\documentclass[11pt]{article}
\usepackage[ansinew]{inputenc}
\usepackage{a4wide}
\usepackage{calc}
\usepackage{times}
\usepackage{fancyhdr}
\usepackage{amssymb,amsmath,amstext,amsthm,amsfonts, ams fonts}
\usepackage{colortbl,color,graphics,graphicx,epsfig}
\definecolor{gray}{rgb}{0.8,0.8,0.8}
\graphicspath{{pics/}}
\usepackage{dsfont}
\usepackage{bbm}
\usepackage{nicefrac}
\usepackage{natbib}
\usepackage[pdflinkmargin=5pt,pdfstartview={FitBH -32768},plainpages=false]{hyperref}
\usepackage{comment,lscape,longtable}

\renewcommand{\P}{\mathbb{P}}
\newcommand{\F}{\mathcal{F}}

\newcommand{\E}{\mathbb{E}}

\newcommand{\R}{\mathds{R}}
\newcommand{\1}{\mathbbm{1}}
\renewcommand{\d}{\Delta}
\newcommand{\KLEINO}{{\scriptstyle{\mathcal{O}}}}

\renewcommand{\vec}{\operatorname{vec}}

\DeclareMathAccent{\verywidehat}{\mathord}{largesymbols}{'144}

\newcommand{\Acov}{\mathbb{A}\mathbb{C}\textnormal{O\hspace*{0.02cm}V}}
\newcommand{\cov}{\mathbb{C}\textnormal{o\hspace*{0.02cm}v}}
\DeclareMathOperator{\AVAR}{\mathbf{AVAR}}

\newtheorem{remark}{Remark}
\newtheorem{defi}{Definition}
\newtheorem{assump}{Assumption}
\newtheorem{theorem}{Theorem}[section]
\newtheorem{prop}{Proposition}[section]

\newtheorem{cor}[prop]{Corollary}
\newtheorem{assumpsec}[prop]{Assumption}

\allowdisplaybreaks[3]
\begin{document}
\global\long\def\inc#1{\E\left[{#1}\Big|\F_{t_{i-1}}\right]}


\title{Inference for multi-dimensional high-frequency data with an application to conditional independence testing}

\fancyhf{}
\lhead[\thepage ~~~~\textsl{M. Bibinger \& P. A. Mykland}]{}
\rhead[]{\textsl{Inference for multi-dimensional high-frequency data}~~~~ \thepage}
\begin{center}
\huge \noindent
\textbf{Inference for Multi-Dimensional High-Frequency Data: Equivalence of Methods, Central Limit Theorems, and
an Application to Conditional Independence Testing}\\[.5cm]
\Large  Markus Bibinger \& Per A. Mykland\\[.25cm]
\large \noindent
\textsl{Humboldt-Universität zu Berlin and Department of Statistics, University of Chicago}
\end{center}
 \normalsize\textbf{ABSTRACT.
We find the asymptotic distribution of the multi-dimensional multi-scale and kernel estimators for high-frequency financial data with microstructure.
Sampling times are allowed to be asynchronous and endogenous. 
In the process, we show that the classes of multi-scale and kernel estimators for smoothing noise perturbation are asymptotically equivalent in the sense of having the same asymptotic distribution for corresponding kernel and weight functions. 
The theory leads to multi-dimensional stable central limit theorems and feasible versions. 
Hence they allow to draw statistical inference for a broad class of multivariate models which paves the way to tests and confidence intervals in risk measurement for arbitrary portfolios composed of high-frequently observed assets. As an application, we enhance the approach to construct a test for investigating hypotheses that correlated assets are independent conditional on a common factor.} \\[.5cm]
\normalsize\noindent


\section{Introduction\label{sec:1}}
The estimation of daily integrated volatility and covolatility\footnote{Also known as integrated variance and covariance, but we here stick to the more heavily used terminology.} has become a key topic of statistics of high-frequency data and a central building block in model calibration for financial risk analysis. Recent years have seen a tremendous increase in trading activities along with ongoing buildup of computer-based trading. The broad availability of recorded asset prices at such high frequencies magnifies the appeal of statistical methods to efficiently exploit information from the high-frequency data.\\
This article contributes to this strand of literature by considering a continuous semimartingale
\begin{align}\label{X}X_t=X_0+\int_0^t\mu_s\,ds+\int_0^t\sigma_s\,dW_s~,t\in\mathds{R}^+\,,\end{align}
with drift $\mu$, volatility $\sigma$ and a standard Brownian motion $W$, comprising current stochastic volatility models, observed over a fixed time span $[0,T]$ on a discrete grid and by investigating asymptotics when the mesh size of the grid tends to zero. The natural estimator for the quadratic variation (integrated volatility) from equidistant observations of $X$ at $iT/n,i=0,\ldots,n$, is the discrete version called the realized volatility. In the one-dimensional framework, it gives a consistent estimator which weakly converges with usual $\sqrt{n}$-rate to a mixed normal distribution where twice the integrated quarticity occurs as random asymptotic variance (cf.\,\cite{bn3}, \cite{jacodprotter}, \cite{zhang2001}).
Therefore, the concept of stable weak convergence by \cite{renyi} has been called into play to pave the way for statistical inference and confidence intervals. In our setting, stable convergence is equivalent to joint weak convergence with every measurable bounded random variable\footnote{For a discussion of the general case, see p. 270 of \cite{jacodprotter}.} and thus, accompanied by a consistent estimator of the asymptotic variance, allows to conclude a feasible central limit theorem. This reasoning makes stable convergence a key element in high-frequency asymptotic statistics.\footnote{See Section 2 for definition and further discussion.}

The aspiration to progress to more complex statistical models in this research area, has been mainly motivated by economic issues. First of all, in a multi-dimensional framework, different assets are usually not traded and recorded at synchronous sampling times, but geared to individual observation schemes. Employing simple interpolation approaches has led to the so-called Epps effect (cf.\,\cite{epps}) that covariance estimates get heavily biased downwards at high frequencies by the distortion from an inadequate treatment of non-synchronicity.
In the absence of microstructure, the estimator by \cite{hy} remedies this flaw of naively interpolated realized covolatilities and a feasible central limit theorem has been attained in \cite{hy3}.
For synchronous equidistant high-frequency observations of \eqref{X} increasing sample sizes are expected to render the estimation error by discretization smaller and smaller. Contrary to the feature of the statistical model, in many situations high-frequency financial data exhibit an exploding realized volatility when the sampling frequency is too high.\footnote{This is usually seen with the help of a so-called signature plot, see \cite{abdl00} and also the discussion in Chapter 2.5.2 of \cite{lamanga12}.}
This effect is ascribed to market microstructure frictions as bid-ask spreads and trading costs. A favored way to capture this influence is to extend the classical semimartingale model, where the semimartingale acts to describe dynamics of the evolution of a latent efficient log-price which is corrupted by an independent additive noise. Following this philosophy from \cite{zhangmykland}, several integrated volatility estimators have been designed which smooth out noise contamination first. The optimal minimax convergence rate for this model declines to $n^{1/4}$, what is known from the mathematical groundwork provided by \cite{gloter}. This rate can be attained using the multi-scale realized volatility by \cite{zhang}, pre-averaging as described in \cite{JLMPV}, the kernel estimator by \cite{bn2} or a Quasi-Maximum-Likelihood approach by \cite{xiu}. Though the estimators have been found in independent works and rely on various principles, it turned out that they are in a certain asymptotic sense equivalent which is clarified in Section \ref{sec:3} below. \\
Recently, methods to deal with noise and non-synchronicity in one go have been established in the literature. In fact, to each of the abovementioned smoothing techniques (at least) one extension to non-synchronous observation schemes has been proposed. First, the multivariate realized kernels by \cite{bn1} using refresh time sampling are eligible to estimate integrated volatility matrices and guarantee for positive semi-definite estimates at the cost of a sub-optimal convergence rate. \cite{sahalia} suggested to combine a generalized synchronization algorithm with the Quasi-Maximum-Likelihood approach. \cite{lintonpark} use Fourier methods on the same problem. Eventually, a feasible asymptotic distribution theory for the general non-synchronous and noisy setup has been provided by \cite{bibinger2} and \cite{PHY} for hybrid approaches built on the Hayashi--Yoshida estimator and the multi-scale and pre-average smoothing, respectively. Although these estimators combine similar ingredients they behave quite differently, since for the approach in \cite{bibinger2} interpolation takes place on the high-frequency scale after smoothing is adjusted with respect to a synchronous approximation whereas \cite{PHY} suggest to denoise each process first and take the Hayashi-Yoshida estimator from pre-averaged blocks which results in interpolation with respect to a lower-frequency scale.\\
The presented limit theorems and asymptotic distributions of the above estimators in the literature are univariate, i.e.\,only the asymptotic variances of (co-)variation estimators are established. An apparent problem pertinent to applications is to quantify the risk of a collection of high-frequently observed assets.
When $X$ in \eqref{X} is $d$-dimensional, for instance, estimating the quadratic variation of some portfolio as $w_1X^{(1)}+w_2X^{(2)}$ with weights $w_1,w_2$ is based on estimates for the integrated volatilities and the integrated covolatility. As the three estimates are correlated, we are in need of a multivariate limit theorem to deduce the asymptotic variance of the compound estimator.
In this work, we establish multivariate stable limit theorems with the asymptotic variance-covariance matrix of the respective estimation methods along with feasible versions.
Thereto, beyond techniques from statistics of high-frequency data, we exploit elements of matrix calculus. Introduce the multivariate notation by the stable central limit theorem for the realized volatility matrix from regular observations as estimator of the integrated volatility matrix $\int_0^T\Sigma_s\,ds,\Sigma=\sigma\sigma^{\top}$:
\begin{align}\label{cltrv}\hspace*{-.2cm}n^{\frac12}\hspace*{-0.05cm}\vec\hspace*{-0.025cm}\left(\sum_{i=1}^n\hspace*{-0.05cm}\big(X_{\frac{iT}{n}}\hspace*{-0.025cm}-\hspace*{-0.025cm}X_{\frac{(i-1)T}{n}}\big)\big(X_{\frac{iT}{n}}\hspace*{-0.025cm}-\hspace*{-0.025cm}X_{\frac{(i-1)T}{n}}\big)^{\top}\hspace*{-.125cm}-\hspace*{-0.1cm}\int_0^T\hspace*{-0.15cm}\Sigma_s\,ds\hspace*{-0.05cm}\right)\hspace*{-0.025cm}\stackrel{st}{\rightarrow} \hspace*{-0.05cm}MN\Big(0,T\hspace*{-.1cm}\int_0^T\hspace*{-.15cm}\big(\Sigma_s\otimes\Sigma_s\big)\mathcal{Z}\,ds\hspace*{-0.05cm}\Big).\end{align}
The $\vec$-operator transforms the $(d\times d)$ matrix on the left-hand side into a $d^2$-dimensional vector by stacking the columns below each other:
\begin{align*}
&\vec(A)=\left(A^{(11)},A^{(21)},\ldots,A^{(d1)},A^{(12)},A^{(22)},\ldots,A^{(d2)},
\ldots,A^{(d(d-1))},A^{(dd)}\right)^{\top}\in{\mathds{R}}^{d^2}\,,
\end{align*}
for $A=\big(A^{(pq)}\big)_{1\le p,q\le d}\in\R^{d\times d}$.
The mixed normal limit right-hand side comprises a $(d^2\times d^2)$ random asymptotic variance-covariance matrix with the Kronecker square of $\Sigma$. The Kronecker product $A\otimes B\in\R^{d^2\times d^2}$ for $A,B\in\R^{d\times d}$ is defined by
\[ (A\otimes B)^{(d(p-1)+q,d(p'-1)+q')}=A^{(pp')}B^{(qq')},\quad p,q,p',q'=1,\ldots,d.\]
The matrix $\mathcal{Z}$ describes the variance-covariance structure of the empirical covariance matrix of a standard Gaussian vector
\begin{align}\label{EqZ}
{\cal Z}=\cov(\vec(ZZ^\top))\in\R^{d^2\times d^2}\text{ for }Z\sim N(0,E_d)\,,
\end{align}
with $E_d$ the $(d\times d)$ identity matrix. $\cal Z$ is explicit, i.e.\,with $\delta_{p,q}=\1_{\{p=q\}}$:
\[ {\cal Z}^{(d(p-1)+q,d(p'-1)+q')}=(1+\delta_{p,q})\delta_{\{p,q\},\{p',q'\}}, p,q,p',q'=1,\ldots,d,
\]
by the property ${\cal Z}\vec(A)=\vec(A+A^\top)$ for all $A\in\R^{d\times d}$. The matrix $\mathcal{Z}$ is twice the so-called symmetrizer matrix from \cite{magnus}. For realized volatilities to estimate $\int_0^T\sigma_s^2\,ds$ with $\sigma_s$ one-dimensional, we recover their well-known asymptotic variance $2T\int_0^T\sigma_s^4\,ds$. In a two-dimensional setup with volatilities $\sigma_s^{(1)},\sigma_s^{(2)}$ and a correlation process $\rho_s$, we derive as limit variance of the realized covolatility $T\int_0^T(1+\rho_s^2)\big(\sigma_s^{(1)}\sigma_s^{(2)}\big)^2\,ds$. Less familiar are the limiting covariances between realized volatility and realized covolatility $2T\int_0^T\rho_s\big(\sigma_s^{(1)}\big)^3\sigma_s^{(2)}\,ds$ and symmetrically. The form of the asymptotic variance-covariance in \eqref{cltrv} is proved in Appendix \ref{sec:9.1} to prepare for the proofs of the main results.\\
We find the multivariate limit theorems for the multivariate multi-scale estimator under noise in Theorem \ref{procovms} and for the generalized multi-scale estimator under noise and non-synchronous observations in Theorem \ref{procovgms}, which by our equivalence result applies in the same way to realized kernels. We show that endogenous observation times do not complicate the asymptotic law and shed light on the different impact of endogeneities in models with and without microstructure, respectively.\\
Relying on the asymptotic distribution of the considered quadratic covariation matrix estimators, we strive to design a statistical test
for investigating hypotheses, if two processes have zero covariation conditioned on a third one. We end up with a feasible stable central limit theorem for the test statistic involving products of estimators and thus obtain an asymptotic distribution free test. This test, which we call conveniently conditional independence test, renders information about the dependence structure in multivariate portfolios and can be applied to test for zero covariation of idiosyncratic factors in typical portfolio dependence structure models, as the one by \cite{eberlein}. In particular, we may identify dependencies between single assets not carried in common macroeconomic factors that influence the whole portfolio and disentangle those from correlations induced by market influences. \\
The outline of the article is as follows. In Section \ref{sec:2} we first unify the asymptotic analysis of quadratic covariation estimation under noise by proving equivalence of methods. Then, multivariate stable limit theorems are developed. Section \ref{sec:3} proceeds to statistical experiments with noise and non-synchronous endogenous observation times. The conditional independence test is introduced in Section \ref{sec:4} and applied in an empirical study in Section \ref{sec:5} to high-frequency financial data. The proofs can be found in the Appendix.

\section{Estimating the Quadratic Covariation Matrix in Presence of Microstructure Noise\label{sec:2}}
\begin{assump}\label{eff}
Consider a continuous $d$-dimensional It\^{o} semimartingale \eqref{X} adapted with respect to a right-continuous and complete filtration $(\mathcal{F}_t)$ on a filtered probability space $(\Omega,\mathcal{F},(\mathcal{F}_t),\P)$ with adapted locally bounded drift process $\mu$, a $d$-dimensional $(\mathcal{F}_t)$-Brownian motion $W$ and adapted $(d\times d^{\prime})$ c\`{a}dl\`{a}g volatility process $\sigma$. Suppose that $\sigma$ itself is a continuous It\^{o} semimartingale again, given by an equation similar to \eqref{X}. The processes $\sigma$ and $W$ can be dependent, allowing for leverage effect.
\end{assump}
\begin{assumpsec}\label{dis3}
The $d$-dimensional continuous semimartingale $X$ from \eqref{X} is discretely observed on $[0,T]$ with additive noise:
$$ Y_{j}=X_{t_j}+\epsilon_{j}~,j=0,\ldots,n~.$$
The synchronous observation times $t_j,0\le j\le n$, satisfy
\begin{equation}\label{eqdis1}\delta_n=\sup_{j}\left(\left(t_j-t_{j-1}\right),t_0,T-t_n\right)=\mathcal{O}\left(n^{-\frac89-\alpha}\right)\end{equation}
for a constant $0<\alpha\le 1/9$, stating that we allow for a maximum time instant tending to zero slower than with $n^{-1}$, but not too slow.
The microstructure noise is given as a discrete-time process for which the observation errors are assumed to be i.\,i.\,d.\,and independent of the efficient process $X$. Furthermore, the errors have mean zero, and eighth moments exist.
\end{assumpsec}
The variance-covariance matrix of $\epsilon_j,0\le j\le n$, is denoted by ${\mathbf{H}}$ and $\cov\big(\epsilon_j\epsilon_j^{\top}\big)={\mathbf{H}}^{\otimes}\mathcal{Z},0\le j\le n$. In case that $\epsilon_j\sim N(0,\mathbf{H})$, we already know that $\mathbf{H}^{\otimes}=(\mathbf{H}\otimes \mathbf{H})$, but we allow for much more general noise. We write
\begin{align}\Delta_j Y=Y_{t_j}-Y_{t_{j-1}}~\mbox{and}~ \Delta_j^i Y=Y_{t_j}-Y_{t_{j-i}},1\le j \le n,2\le i \le j\,,
\end{align}
for the increments and for increments to longer lags, respectively. Since notation varies between papers, note the correspondence to the other main form:
\begin{align*} \Delta_j Y \mbox{ is the same as } \Delta Y_{t_j}\,.
\end{align*}
An i.\,i.\,d.\,assumption on the noise is standard in related literature, an extension to $m$-dependence and mixing errors can be attained as in \cite{zhangmykland2}. For notational convenience and to find the multivariate analogues of known one-dimensional asymptotic variances of considered estimators, we also restrict ourselves to i.\,i.\,d.\,noise here. Increments in this microstructure noise model
$$\d_j Y=\int_{t_{j-1}}^{t_j}\mu_s\,ds+\int_{t_{j-1}}^{t_j}\sigma_s\,dW_s+\epsilon_j-\epsilon_{j-1}$$
are substantially governed by the noise, since any component of the second addend is $\mathcal{O}_{\P}(\delta_n^{1/2})$ and the drift acts only as nuisance term of order in probability  $\mathcal{O}_{\P}(\delta_n)$ for each component. For an accurate estimation of the quadratic covariation matrix in the presence of noise smoothing methods are applied. We now discuss several main approaches and integrate them in a unifying theory. To this end, we show that two prominent methods are asymptotically equivalent.\\
The asymptotic distributions of considered estimators hinge on the random volatility process $\sigma_s$. Thus, stable weak convergence is an essential concept.\footnote{Let $Z_n$
be a sequence of ${\mathcal X}$-measurable random variables, with $\F_T \subseteq \chi$.
We say that $Z_n$ {\it converges stably in law to} $Z$ as $n \rightarrow \infty$
if $Z$ is measurable with respect to an extension
of ${\mathcal X}$ so that for all $A \in \F_T$ and for all
bounded continuous $g$, $EI_A g(Z_n) \rightarrow EI_A g(Z)$ as
$n \rightarrow \infty$. $I_A$ denotes the indicator function of $A$, and $=1$ if $A$ and $=0$ otherwise. In the case of no
microstructure, ${\mathcal X} = \F_T$. If there is microstructure, ${\mathcal X}$ is formed as the smallest sigma-field
containing $\F_T$ and also making the microstructure measurable.
We refer to \cite{jacod1} and \cite{jacodprotter} for background information on stable convergence for this estimation problem.}
Stable central limit theorems allow for feasible limit theorems and hence confidence if the asymptotic variance-covariance matrix can be estimated consistently.
\footnote{Stable convergence also permits the suppression of drift through
measure change,  see
Section 2.2 of \cite{myklandzhang_Mfinite}, which draws on \cite{rootzen80}. The device is similar to
the passage to risk neutral measures in finance, going back to \cite{ross76} and \cite{harrisonkreps79}. 
This mode of convergence also permits
the localization of processes such as volatility, so they can be assumed bounded, see Chapter 2.4.5 of \cite{lamanga12}.}

\subsection{The Multivariate Multi-Scale and Kernel Estimators}
For the estimation of the quadratic variation the following rate-optimal estimators with similar asymptotic behavior have been proposed in the literature: the multi-scale approach by \cite{zhang}, pre-averaging by \cite{JLMPV}, the kernel estimator by \cite{bn2} and a Quasi-Maximum-Likelihood estimator by \cite{xiu}. We investigate the variance-covariance structure of the multivariate multi-scale estimator explicitly, but since all these estimators have a similar structure as quadratic form of the discrete observations, analogous reasoning will apply to the other methods. In particular, we shed light on the connection to the kernel approach to profit at the same time from the considerations by \cite{bn2} pertaining parametric efficiency and the asymptotic features of different kernel functions.
The multivariate multi-scale estimator
\begin{align}\label{ms}\widehat{\left[ X, X\right]}_T^{(multi)}=\sum_{i=1}^{M_n}\frac{\alpha_{i}}{i}\sum_{j=i}^n\d_j^iY(\d_j^iY)^{\top}\end{align}
arises as linear combination of averaged lower-frequent realized volatility matrices using frequencies $i=1,\ldots,M_n$. Estimator \eqref{ms} is the multi-dimensional version of the estimator from \cite{zhang}.\\
For discrete weights $\alpha_i,1\le i\le M_n$, with $\sum_{i=1}^{M_n}\alpha_i=1$ and \(\sum_{i=1}^{M_n}(\alpha_i/i)=0\), the expression
\begin{align}\label{weights}
\alpha_{i}&=\frac{i}{M_n^{2}}h\left(\frac{i}{M_n}\right)-\frac{i}{2M_n^{3}}h^{\prime}\left(\frac{i}{M_n}\right)+\frac{i}{6M_n^4}(h^{\prime}(1)-h^{\prime}(0))-\frac{i}{24M_n^5}(h^{\prime\prime}(1)-h^{\prime\prime}(0))\,,\end{align}
adopted from \cite{zhang}, with twice continuously differentiable functions $h$ satisfying $\int_0^1 xh(x)\,dx$ $=1$ and $\int_0^1 h(x)\,dx=0$, gives access to a tractable class of estimators. The multi-scale frequency is chosen $M_n=c\,\sqrt{n}$ with a constant $c$, minimizing the overall mean square error to order $n^{-1/4}$. The estimator is thus rate-optimal according to the lower bounds for convergence rates by \cite{gloter} and \cite{bibinger}.\\ At the present day, it is commonly known that the nonparametric smoothing approaches to cope with noise contamination have a connatural structure and related asymptotic distributions. A prominent intensively studied alternative to the multi-scale approach is the (realized) kernel estimator
\begin{align}\label{kernel}\notag\widehat{\left[ X, X\right]}_T^{(kernel)}&=\sum_{j=1}^n\d_j Y(\d_j Y)^{\top}\\
&\quad +\sum_{h=1}^{H_n}\mathfrak{K}\left(\frac{h}{H_n}\right)\Big(\sum_{j=h+1}^n \d_jY(\d_{j-h}Y)^{\top}+\d_{j-h}Y(\d_jY)^{\top}\Big)\,,\end{align}
with a four times continuously differentiable kernel $\mathfrak{K}$ on $[0,1]$, which satisfies the following conditions:
\[ \hspace*{-0.2cm}\max{\left\{\int_0^1 \hspace*{-0.1cm}\mathfrak{K}^2(x)\,dx,\int_0^1  \hspace*{-0.1cm}(\mathfrak{K}^{\prime}(x))^2\,dx,\int_0^1 \hspace*{-0.1cm} (\mathfrak{K}^{\prime\prime}(x))^2\,dx\right\}} \hspace*{-0.1cm}<\infty, \mathfrak{K}(0)=1, \mathfrak{K}(1)=\mathfrak{K}^{\prime}(0)=\mathfrak{K}^{\prime}(1)=0.\]
This is the multi-dimensional version of the estimator by \cite{bn2}. In the one-dimensional setup \eqref{kernel} has been motivated as linear combination of realized autocovariances of the discretely observed process.\\ The subsequent explicit relation between kernel and multi-scale estimator enables us to embed the findings about several kernels and the construction of an asymptotically efficient one for the parametric model provided by \cite{bn2}. Since the multi-scale approach exhibits good finite-sample properties in the treatment of end-effects, it can be worth to road-test resulting transferred multi-scale estimators in practice.

\subsection{Asymptotic Equivalence of the Multi-Scale and Kernel Estimators}
The multi-scale and kernel estimators defined in \eqref{ms} and \eqref{kernel} are sensitive to end-effects which is caused by the dominating noise component whose variance-covariance matrix ${\mathbf{H}}$ does not depend on $n$. Due to end-effects, on Assumption \ref{dis3}, the estimators \eqref{ms} and \eqref{kernel} with weights determined by \eqref{weights} and corresponding kernels have a bias $-2\,{\mathbf{H}}$ and $2\,{\mathbf{H}}$, respectively. We here investigate a correction to each of the two types of estimator:\\
{\it Correction to Multi-scale:} Follow \cite{zhang} by modifying the first two weights
\begin{equation}
\alpha_1\mapsto \alpha_1+2/n,\alpha_2\mapsto \alpha_2-2/n, (\alpha_i)_{3\le i\le M_n}\mapsto (\alpha_i)_{3\le i\le M_n} .
\end{equation}
{\it Correction to the Kernel estimator:}
\begin{equation}
\mbox{ multiplying the realized volatility matrix in the first addend with } \frac{n-1}{n} .
\end{equation}
This correction is different from the `jittering' approach provided in \cite{bn2}.\footnote{Section 2.6 p. 1487-88 of  \cite{bn2}.} The bias-corrections do not affect the asymptotic variance-covariance structure of the estimators.
We call the adjusted estimators, respectively,
\begin{equation*}
\widehat{\left[ X, X\right]}_T^{(multi,adj)} ~\mbox{ and  }~~
\widehat{\left[ X, X\right]}_T^{(kernel,adj)} .
\end{equation*}
We then obtain the following direct asymptotic equivalence of the two estimators.
\begin{theorem}\label{thm:comparison-raw}
For each kernel function $\mathfrak{K}$ matching the assertions above, for the estimators defined in \eqref{ms} and \eqref{kernel} with weights determined by \eqref{weights} and $h=\mathfrak{K}^{\prime\prime}$, we have
\begin{align} n^{\frac{1}{4}}\left(\widehat{\left[ X, X\right]}_T^{(multi)}-\widehat{\left[ X, X\right]}_T^{(kernel)}+4\,{\mathbf{H}}\right)\stackrel{p}{\rightarrow}0\,,\end{align}
as $n\rightarrow\infty$, $M_n=H_n\rightarrow\infty$. The term $4\,{\mathbf{H}}$ is only due to the different impact of end-effects.
For the bias-corrected versions, we thus conclude that
\begin{align} n^{\frac{1}{4}}\left(\widehat{\left[ X, X\right]}_T^{(multi,adj)}-\widehat{\left[ X, X\right]}_T^{(kernel,adj)}\right)\stackrel{p}{\rightarrow}0\,,\end{align}
as $n\rightarrow\infty$, $M_n=c_{multi}\sqrt{n}$ and $H_n=c_{kern}\sqrt{n}$.
\end{theorem}
\begin{remark} (Dependent noise.)
\label{rmk:m-dep}
In the case of $m$-dependence it will be convenient to discard the first $m$ frequencies and renormalize in \eqref{ms}. The adjusted estimator is robust.
\end{remark}
\begin{table}[t]\renewcommand{\arraystretch}{1.5}
\begin{center}
\begin{tabular}{|c|c|}
\hline
\rowcolor{gray}kernel & $\mathfrak{K}$ \\
\hline
cubic &	$1-3x^2+2x^3$\\
Parzen &	$(1-6x^2+6x^3)\1_{\{x\le 1/2\}}+2(1-x)^3\1_{\{x>1/2\}}$	\\
$r$th Tukey-Hanning &$\sin\left(\frac{\pi}{2}(1-x)^r\right)^2$\\
\hline
\rowcolor{gray}kernel & first-order weights ~$\alpha_i$\\
\hline
cubic &	$\frac{12i^2}{(M)^3}-\frac{6i}{(M)^2}$\\
Parzen & $\frac{i}{M^2}\left(\frac{36i}{M}-12\right)$ for $i\le M/2$ and $\frac{i}{M^2}\left(12-\frac{12i}{M}\right)$ for $i>M/2$\\
$r$th Tukey-Hanning &	$\frac{\pi i r (1-\frac{i}{M})^{r-2} \left((r-1) \sin{(\pi (1-\frac{i}{M})^r)}+\pi r(\frac{i}{M}-1)^r \cos{(\pi (1-\frac{i}{M})^r)}\right)}{2 M^2}$\\
\hline
\end{tabular}\end{center}
\caption{\label{tab1}Collection of important kernels and corresponding weights for the multi-scale (first order term).}
\end{table}
\begin{remark} (Strong representation.)
The result of Theorem \ref{thm:comparison-raw} is similar to other ``strong representation'' results in the high-frequency
literature,
such as in \cite{lancov2011} (see key equation (39) on p. 41) and \cite{efficientbipower}, Theorem 4. (The convergence is in probability,
but is comparable to strong representation through a standard subsequence-of-subsequence argument.)
\end{remark}
Since the motivation of the multi-scale and the kernel approach is quite different, the asymptotic equivalence in Theorem \ref{thm:comparison-raw} is an intriguing result. The equivalence and its proof also reveal how refinements and results for one estimator can be transferred to the other.

\subsection{Optimal Choice of Weights, and Asymptotic Distribution}
The standard weights employed in \cite{zhang}
\begin{align}\label{cubweights}
\alpha_{i}&=\left(\frac{12i^2}{(M_n^3-M_n)}-\frac{6i}{(M_n^2-1)}-\frac{6i}{(M_n^3-M_n)}\right)=\frac{12i^2}{M_n^3}-\frac{6i}{M_n^2}\left(1+\KLEINO(1)\right)\end{align}
minimize the variance by noise and lead to, as mentioned by \cite{bn2}, the same asymptotic properties as for the kernel estimator \eqref{kernel} with a cubic kernel. However, as derived by \cite{bn2} there are kernels surpassing the cubic kernel in efficiency by shrinking the signal and cross parts of the variance while allowing for an increase in the noise variance and striving for the best balance of all three. A fourth term appearing in the asymptotic (co-)variances, see \eqref{covms} below, induced by end-effects and noise, can be circumvented by their `jittering'  technique.
Asymptotically, Tukey-Hanning kernels as listed in Table \ref{tab1} combined with this `jittering' can attain the optimal asymptotic variance in the one-dimensional parametric case known from the inverse Fisher information in \cite{gloter}. All weights \eqref{weights} satisfy the relations $\sum_{i=1}^{M_n}\alpha_i=1$ and $\sum_{i=1}^{M_n}\alpha_i/i=0$.
Classical pre-averaging is asymptotically equivalent to the Parzen kernel. This linkage has been shown by \cite{kinnepoldi}; see also the discussion in \cite{JLMPV} (Remark 1, p. 2255). At this stage, we derive the multivariate stable central limit theorem along with the asymptotic variance-covariance matrix for the equidistant observations setup which will be extended to irregular sampling below within our general non-synchronous model.
\begin{table}[t]\renewcommand{\arraystretch}{1.5}
\begin{center}
\begin{tabular}{|c|c|c|c|c|}
\hline
\rowcolor{gray}kernel & $\mathfrak{N}_1^{\alpha}$ & $\mathfrak{D}^{\alpha}$ &$\mathfrak{M}^{\alpha}$ &$\mathfrak{N}_2^{\alpha}$\\
\hline
cubic &	$12$ & $13/70$ &$6/5$ &$6/5$ \\
Parzen &	$24$ & $3/4$ & $151/560$ & $15/40$	\\
$1$st Tukey-Hanning & $\pi^4/8$ &$\pi^2/16$&$3/8$& $\pi^2/8$\\
$16$th Tukey-Hanning & $14374$ &$5.132$&$0.0317$& $10.264$\\
\hline
\end{tabular}\end{center}
\caption{\label{tab2}Constants in asymptotic covariance for important kernels.}
\end{table}

\begin{theorem}\label{procovms}
On the Assumptions \ref{eff} and \ref{dis3} with $t_i=iT/n,0\le i\le n$, the multi-scale estimator \eqref{ms} with $M_n=c\,\sqrt{n}$, and weights \eqref{weights}, and by the equivalence also the corresponding kernel estimator, obey multivariate stable central limit theorems
\begin{align}n^{\frac14}\,\vec\Big(\widehat{\left[ X, X\right]}_T^{(multi)}-\int_0^T\Sigma_s\,ds\Big)\stackrel{st}{\rightarrow} MN\big(0,\Acov\big)\,,\end{align}
with mixed normal limit distribution and with the asymptotic variance-covariance matrix
\begin{align} \label{covms}\Acov&=4\mathfrak{D}^{\alpha}\,c\,T\int_0^T (\Sigma_s\otimes \Sigma_s)\mathcal{Z}\,ds+2\,\mathfrak{N}_1^{\alpha} \,c^{-3}\left({\mathbf{H}}\otimes{\mathbf{H}}\right)\mathcal{Z}\\
&\notag\quad  +2\,c^{-1}\mathfrak{M}^{\alpha}\,\int_0^T\left({\mathbf{H}}\otimes\Sigma_s+\Sigma_s\otimes{\mathbf{H}}\right)\mathcal{Z}\,ds +2\,c^{-1}\mathfrak{N}_2^{\alpha}\mathbf{H}^{\otimes}\mathcal{Z}\,,\end{align}
with constants $\mathfrak{D}^{\alpha}$, $\mathfrak{N}_1^{\alpha}$, $\mathfrak{N}_2^{\alpha}$ and $\mathfrak{M}^{\alpha}$ depending on the specific kernel, see Table \ref{tab2}.
\end{theorem}
A generalization for non-equidistant sampling is covered by Theorem \ref{procovgms} in Section \ref{sec:3}. In this case, the first addend of \eqref{covms} (signal term) hinges on a function \eqref{wlasa}, while the other terms are analogous.
In the special case $d=1$, we obtain the asymptotic variance of the one-dimensional multi-scale estimator as given in \cite{zhang}. The last addend involving $\mathbf{H}^{\otimes}\mathcal{Z}=\cov\big(\epsilon_j\epsilon_j^{\top})$ is induced by end-effects and noise and can be circumvented by the jittering technique, see Section 2.6 of \cite{bn2}. For the cross terms note the identity
\[\left({\mathbf{H}}\otimes\Sigma_s+\Sigma_s\otimes{\mathbf{H}}\right)\mathcal{Z}=\mathcal{Z}\left({\mathbf{H}}\otimes\Sigma_s\right)\mathcal{Z}=\mathcal{Z}\left(\Sigma_s\otimes {\mathbf{H}}\right)\mathcal{Z}\,.\]

\section{Estimating the Integrated Volatility Matrix under Asynchronicity and Microstructure\label{sec:3}}
\subsection{Asymptotic Distribution of the Generalized Multi-Scale Estimator}
This section focuses on the general model -- comprising non-synchronous observation times and noise perturbation -- and an hybrid approach founded on a combination of the estimators from Section \ref{sec:2} and the estimator for non-synchronous non-noisy observations by \cite{hy}. First, observation times are deterministic or random and independent of $Y$. In Section \ref{sec:3.1} robustness against endogenous sampling is established.
\begin{assumpsec}\label{dis4}
The process $X$ is observed non-synchronously with additive microstructure noise:
$$Y^{(p)}_{t_j^{(p)}}=X^{(p)}_{t_j^{(p)}}+\epsilon^{(p)}_{j}~,j=0,\ldots,n_p,p=1,\ldots,d,~\text{on}~[0,T]\,.$$
The sequences of observation times are regular in the sense that $n_p/n_q\rightarrow K_{pq}$ with constants $0<K_{pq}<\infty$. For a constant $0<\alpha\le 1/9$, it holds that
\begin{align}\label{eqdis3}
\delta_n&=\sup_{(i,p)}{\left(\left(t_i^{(p)}-t_{i-1}^{(p)}\right),t_0^{(p)},T-t_{n_p}^{(p)}\right)}~\;\,=\mathcal{O}\left(\sup_p(n_p)^{-\frac{8}{9}-\alpha}\right)~.\end{align}
The observation errors are i.\,i.\,d.\,sequences, independent of the efficient processes, centered and eighth moments exist. Noise components can be mutually correlated only at synchronous observations.
\end{assumpsec}
We work conditional given the exogenous observation times. In the following we establish the asymptotic variance-covariance matrix for a generalized multi-scale method proposed in \cite{bibinger} and \cite{bibinger2}. It arises as a convenient composition of the multi-scale estimator from Section \ref{sec:2} and a synchronization approach inspired by the estimator suggested in \cite{hy}.\\
To handle non-synchronicity, introduce the notion of next- and previous-tick interpolations:
\begin{align*}t_p^+(s)=\min_{i\in\{0,\ldots,n_p\}}{\left(t_i^{(p)}|t_i^{(p)}\ge s\right)}~\mbox{and}~t_p^-(s)=\max_{i\in\{0,\ldots,n_p\}}{\left(t_{i}^{(p)}|t_i^{(p)}\le s\right)}\end{align*}
for $p=1,\ldots,d$, and $s\in[0,T]$. An important synchronous grid is given by the refresh times introduced in \citet{bn1}:
\begin{align*}T_0=\max_p{\left(t_p^+(0)\right)}\,,\,T_i=T_{i-1}+\max_p{\left(t_p^+(T_{i-1})\right)}\,,i=1,\ldots,N\,.\end{align*}
For the construction of the estimator, virtually we can think of an idealized synchronous approximation given by the $(N+1)$ refresh times, apply subsampling and the multi-scale extension to this scheme, and afterwards interpolate to the next observed values on the highest available frequency. This generalized multi-scale estimator is
\begin{align}\label{gms}\widehat{\left[ X, X\right]}_T^{(multi)}&=\sum_{i=1}^{M_N}\frac{\alpha_i}{i}\sum_{j=i}^{N}\left(Y_{T_j}^+-Y_{T_{j-i}}^-\right)\left(Y_{T_j}^+-Y_{T_{j-i}}^-\right)^{\top}\,,\\
\nonumber \mbox{with}~~Y_{T_j}^+&=\Big(Y_{t_p^+(T_j)}^{(p)}\Big)^{\top}_{1\le p\le d}\,,\, Y_{T_j}^{-}=\Big(Y_{t_p^-(T_j)}^{(p)}\Big)^{\top}_{1\le p\le d}\,,j=0,\ldots,N.\end{align}
Without loss of generality suppose all next-tick interpolations for $j=N$ and previous-tick interpolations for $j=0$ exist (simply exclude the first and last refresh time else).
This estimator crucially differs from the approach by \cite{PHY}, which mimics the form of the traditional Hayashi-Yoshida estimator, but bound to a low-frequency scheme of pre-averaged observations over blocks of order $\sqrt{n}$ high-frequency observations. The estimator \eqref{gms} relies more on the principle of the refresh-time approximation and exhibits a simpler and for most setups much smaller variance. Contrarily to \cite{bn1}, we utilize pre- and next-tick interpolations such that the final estimator has no bias due to non-synchronicity. For the reason of various estimators in the general model using different compositions of the methods, the article on hand can not accomplish a unified theory that is applicable to all alternative approaches as \cite{sahalia}, \cite{bn1} and \cite{PHY}. Unlike their roots from Section \ref{sec:2} they are not asymptotically equivalent any more. We consider \eqref{gms} because the method attains a much smaller discretization variance in comparison to the one by \cite{PHY}, is rate-optimal and a feasible univariate central limit theorem is accessible from \cite{bibinger2}.
\begin{remark} (Identical results for kernel estimators.)
Since equations \eqref{eqdis1} and \eqref{eqdis3} are the same, it follows from Section \ref{sec:2} that our results on irregular sampling
for the synchronous case, where the generalized multi-scale estimator \eqref{gms} coincides with the original one \eqref{ms}, in the following apply identically to kernel estimators. Furthermore, all results for the estimator \eqref{gms} apply to a generalized kernel estimator with refresh time sampling as in \eqref{gms}.
\end{remark}

\begin{defi}\label{lsadefi}For observation times $t_j^{(p)},0\le j\le n_p,1\le p\le d$, define the functional sequences
\begin{align}\label{lsa}\mathfrak{G}_{N,r}(t)=\frac{N}{r}\sum_{T_l\le t}\left(T_l-T_{l-1}\right)\sum_{q=1}^{r\wedge l}\left(T_{l-q+1}-T_{l-q}\right)\,,\end{align}
and $S^N(t)\in\mathds{R}^{d\times d}$ for each $t\in[0,T]$ with entries
\begin{align}\label{dr}\big(S^N(t)\big)^{(pq)}=\frac{1}{N}\sum_{T_l\le t}\Big(\1_{\{t_p^+(T_l)=t_{q}^+(T_l)\}}+2\,\1_{\{t_p^+(T_l)=t_{q}^+(T_{l-1})\}}+\sum_{u=0}^l\1_{\{t_p^+(T_l)=t_{q}^-(T_u)\}}\Big)~.\end{align}
\end{defi}
\begin{assumpsec}\label{aqvt2}
Assume that the sequence $\mathfrak{G}_{N,r}$ from \eqref{lsa} and the sequences \eqref{dr} satisfy the convergences
\begin{enumerate}
\item[(i)]As $N\rightarrow\infty$ and $r\rightarrow\infty$ with $r=\KLEINO(N)$: $\mathfrak{G}_{N,r}(t)\rightarrow \mathfrak{G}(t)$ and $S^N(t)\rightarrow S(t)$, for continuous differentiable (in $t$) limiting functions $\mathfrak{G}$ and $S$ on $[0,T]$.
\item[(ii)]For any null sequence $(h_N),\,h_N=\mathcal{O}\left(N^{-1}\right)$:
\begin{align}\label{lasa}\frac{\mathfrak{G}_{N,r}(t+h_N)-\mathfrak{G}_{N,r}(t)}{h_N}\rightarrow \mathfrak{G}^{\prime}(t)\,,\end{align}
\begin{align}\label{dra}\frac{S^N(t+h_N)-S^N(t)}{h_N}\rightarrow S^{\prime}(t)\,,\end{align}
uniformly on [0,T] as $N\rightarrow \infty$.
\item[(iii)]Assume that for all $p,p',q,q'\in\{1,\ldots,d\}$, the following limits exist
\begin{align}\label{chi}\chi^{pp'}_{qq'}&=\lim_{N\rightarrow\infty}\frac{M_N^3}{N}\sum_{i=1}^{M_N}\frac{\alpha_i^2}{i^2}\sum_{j=i+1}^N\Big(\1_{\{t_p^+(T_j)=t_{p'}^+(T_j)\}}\1_{\{t_q^-(T_{j-i})=t_{q'}^-(T_{j-i})\}}\\
&\nonumber\hspace*{4cm} +2\,\1_{\{t_p^+(T_j)=t_{p'}^+(T_{j-1})\}}\1_{\{t_q^-(T_{j-i-1})=t_{q'}^-(T_{j-i-1})\}}\Big)\,.\end{align}
\item[(iv)]Assume the existence of
\begin{align}\hspace*{-.65cm}\lim_{N\rightarrow\infty}M_N^{-1}\Big(\sum_{j=1}^{M_N}\Big(\1_{\{t_p^+(T_j)=t_{q}^+(T_{j})\}}+2\1_{\{t_p^+(T_j)=t_{q}^+(T_{j-1})\}}\Big) +\sum_{j=N-M_N}^{N}\1_{\{t_p^-(T_j)=t_{q}^-(T_{j})\}}\Big)\,.\end{align}
\end{enumerate}
\end{assumpsec}
The existence of the limit $\mathfrak{G}$ of $\mathfrak{G}_{N,r}$ is essential to establish an asymptotic distribution theory, since it dominates the terms that appear in the (co-)variances of the multi-scale and related estimators and contribute to the asymptotic (co-)variance, namely the following existing limit:
\begin{align}\label{wlasa}\hspace*{-.25cm}D^{\alpha}(t)=\lim_{N\rightarrow\infty}\hspace*{-.05cm}{\left(\frac{N}{M_N}\sum_{T_l\le t}(T_l-T_{l-1})\sum_{i,k=1}^{M_N}\alpha_i\alpha_k\hspace*{-.15cm}\sum_{q=1}^{\min{(l,i,k)}}\hspace*{-.15cm}\left(1-\frac{q}{i}\right)\left(1-\frac{q}{k}\right)(T_{l-q+1}-T_{l-q})\right)}.\end{align}
In the equidistant synchronous setup $D^{\alpha}(t) = \mathfrak{D}^{\alpha} t\,T$, with the constant $\mathfrak{D}^{\alpha}$ found in Theorem \ref{procovms}.

\begin{theorem}\label{procovgms}
On the Assumptions \ref{eff}, \ref{dis4} and \ref{aqvt2}, the generalized multi-scale estimator \eqref{gms} with $M_{N}=c\sqrt{N}$ and weights \eqref{weights} obeys the multivariate stable central limit theorem:
\begin{align}\label{gmsclt}N^{1/4}\left(\widehat{\left[ X, X\right]}_T^{(multi)}-\int_0^T\Sigma_s\,ds\right)\stackrel{st}{\rightarrow}MN\big(0,\Acov\big)\,,\end{align}
with mixed normal limit distribution and with the asymptotic variance-covariance matrix
\begin{align}\label{covgms}\Acov=&4\,c\int_0^T(D^{\alpha})^{\prime}(s)(\Sigma_s\otimes \Sigma_s)\mathcal{Z}\,ds+2\,c^{-3}\left(\mathbf{H}\otimes \mathbf{H}\right)^*\mathcal{Z}\\ &\quad\nonumber
+c^{-1}\,\mathfrak{M}^{\alpha}\int_0^T\left(\tilde{\mathbf{H}}_s\otimes\Sigma_s+\Sigma_s\otimes \tilde{\mathbf{H}}_s\right)\mathcal{Z}\,ds+c^{-1}\, \mathfrak{N}_2^{\alpha}\,\tilde{\mathbf{H}}^{\otimes}\mathcal{Z}\,,\end{align}
with \eqref{wlasa} and the following existing limits:
\begin{subequations}
\begin{align}\label{tildelimit}\big(\tilde{\mathbf{H}}_s\big)^{(pq)}&=\mathbf{H}^{(pq)}\big(S^{\prime}(s)\big)^{(pq)}\,,\\
\label{starlimit}\left(\left(\mathbf{H}\otimes \mathbf{H}\right)^*\right)^{(d(p-1)+q\,,\,d(p'-1)+q')}&=\mathbf{H}^{(pp')}\mathbf{H}^{(qq')}\chi^{pp'}_{qq'}\,,\end{align}
\vspace*{-1.5cm}

\begin{align}
\label{2starlimit}\big(\tilde{\mathbf{H}}^{\otimes}\big)^{(pq)}=\left(\mathbf{H}^{\otimes}\right)^{(pq)}&\Big(\lim_{N\rightarrow\infty}M_N^{-1}\sum_{j=1}^{M_N}\Big(\1_{\{t_p^+(T_j)=t_{q}^+(T_{j})\}}+2\1_{\{t_p^+(T_j)=t_{q}^+(T_{j-1})\}}\Big)\\ &\quad +\lim_{N\rightarrow\infty}M_N^{-1}\sum_{j=N-M_N}^{N}\1_{\{t_p^-(T_j)=t_{q}^-(T_{j})\}}\Big)\nonumber\,,\end{align}
\end{subequations}
for $p,p',q,q'\in\{1,\ldots,d\}$ with $S^{\prime}$ from \eqref{dra} and $\chi^{pp'}_{qq'}$ from \eqref{chi}.
\end{theorem}
In a synchronous setting  $\big(\mathbf{H}\otimes \mathbf{H}\big)^*=\mathfrak{N}^{\alpha}_1\big(\mathbf{H}\otimes \mathbf{H}\big)$ ($\mathfrak{N}^{\alpha}_1=12$ for the cubic kernel), $\tilde{\mathbf{H}}_s=2{\mathbf{H}}$  and $\tilde{\mathbf{H}}^{\otimes}=2\mathbf{H}^{\otimes}$, and then \eqref{covgms} coincides with \eqref{covms} except for the influence of irregular sampling. In particular, the asymptotic variance-covariance matrix of the multi-scale estimator for synchronous but non-equidistant sampling coincides with \eqref{covms}, but in the discretization part the derivative of \eqref{wlasa} analogously defined for the one observation scheme replaces the constant $\mathfrak{D}^{\alpha} T$.\\
Interestingly, in most situations non-diagonal entries of $S(t)$ equal zero as well as $\chi^{pp'}_{qq'}$ whenever $p\ne p'$ or $q\ne q'$, such that the noise part of covariances vanishes. We obtain the following important result for the completely non-synchronous case.
\begin{cor}\label{corcovgms}
In the case that no synchronous observations take place: $t_i^{(p)}\ne t_j^{(q)}$ for all $i,j$ and $p\ne q$ (or the amount of synchronous observations tends to zero as $N\rightarrow\infty$), \eqref{gmsclt} holds and  \eqref{tildelimit}, \eqref{starlimit} and \eqref{2starlimit} simplify to
\begin{subequations}
\begin{align}\label{tildelimitnew}\big(\tilde{\mathbf{H}}_s\big)^{(pq)}&=2\,\mathbf{H}^{(pq)}\delta_{p,q}\\
\label{starlimitnew}\left(\left(\mathbf{H}\otimes \mathbf{H}\right)^*\right)^{(d(p-1)+q\,,\,d(p'-1)+q')}&=\delta_{p,p'}\delta_{q,q'}\mathbf{H}^{(pp)}\mathbf{H}^{(qq)}\mathfrak{N}_1^{\alpha}\\
\big(\tilde{\mathbf{H}}^{\otimes}\big)^{(pq)}&=2\,\big({\mathbf{H}}^{\otimes}\big)^{(pq)}\delta_{p,q}\label{2starlimitnew}\end{align}
\end{subequations}
where $\delta_{p,q}=\1_{\{p=q\}}$.
\end{cor}
\begin{remark}
Our major focus is not on the theoretical limits $\mathfrak{G}$ and of other sequences, since in the general case they
are specified only as limits. We do not need these values, however, for inference, as we shall see in the Section 3.3 on feasible inference.\\
Note that convergence of \eqref{lsa} is the natural assumption to derive a central limit theorem for irregularly spaced (non-equidistant)  observations already in the one-dimensional framework. It emulates the asymptotic quadratic variation of time for realized volatility to an asymptotic long-run variation of time emerging in the variance for subsampling and the other smoothing approaches. Not directly the limit of \eqref{lsa} will appear in the asymptotic variance, but some limiting function additionally involving specific weights (the kernel). If we think of random sampling independent of $Y$, the structure of \eqref{lsa} will be particularly simple for i.\,i.\,d.\,time instants. Virtually, only the expectation will matter and we can apply the standard law of large numbers. Assuming \eqref{lasa} is less restrictive than the assertion in \cite{zhang}, i.\,e.\,sampling needs not to be close to an equidistant scheme in the sense that asymptotic quadratic variation of time converges to T at T. Remarkably, for the popular model of homogenous Poisson sampling independent of $Y$ with expected time instants $T/n$, the asymptotic variance of the integrated volatility estimator is the same as for equidistant observations. This emanates from the i.\,i.\,d.\,nature of time instants and the vanishing influence of the first addend $2T/(nr)$ in \eqref{lsa} as $r\rightarrow \infty$. The finite sample correction factor in \eqref{lasa} for this Poisson setup is thus $(r+1)/r$.
\end{remark}
\begin{remark}(Pairwise refresh times) Instead of subsampling geared to the refresh time scheme in \eqref{gms} we can as well use pairwise refresh times to estimate each entry of the integrated volatility matrix, i.e.\,to estimate $\int_0^T\Sigma_s^{(pq)}\,ds$ we work with refresh times build from $(t_i^{(p)})_{1\le p\le n_p},(t_i^{(q)})_{1\le q\le n_q}$. Especially in case of very different liquidities the pairwise estimation can be more efficient in finite samples. The variance-covariance structure for a pairwise generalized multi-scale estimator is slightly more cumbersome -- but of the same nature as \eqref{covgms}. The discretization variance-covariance terms are analogous where $D^{\alpha}(s)$ is defined as \eqref{wlasa} but based on refresh times built from all involved components (which means from 4 schemes for covariance non-diagonal entries and one or two schemes for diagonal variance entries). In the other parts the covariance entries vanish in the typical setup without synchronous observations. Connatural terms as above show up when we replace in each entry refresh times of all components by refresh times of involved components.
\end{remark}
At first glance the simple appearance of the variance-covariance of generalized multi-scale estimates in the typical general setup where all observations are non-synchronous and in the presence of microstructure noise is intriguing. It hinges only on the discretization error as if we had synchronous observations at the refresh times $T_i,i=0,\ldots,N$. The noise falls out of the asymptotic covariances on the assumption that observation errors at different observation times are independent.\\
This constitutes another nice property of the generalized multi-scale method that a multivariate limit theorem \eqref{gmsclt} is available and covariances are pretty simple. Here, we benefit from the construction of \eqref{gms} where interpolation effects and hence the discretization error due to non-synchronicity is asymptotically negligible. This is in line with the result of \cite{BRHM} that in this general model with microstructure noise and non-synchronicity the noise prevails such that the discretization variance-covariance is asymptotically not affected by non-synchronicity.

\subsection{Incorporating Endogenous Sampling\label{sec:3.1}}
One crucial limitation of the observation model with Assumption \ref{dis4} is that observation times are supposed to be exogenous and not dependent on the process $Y$. This appears unrealistic when observations come at random trading times. A prominent contribution in which volatility estimation in presence of endogenous random observation times has been considered is \cite{fukasawa}, other works dedicated to endogenous sampling include \cite{limyk14} and \cite{hjy}. Especially the limit theorem for realized volatility by \cite{fukasawa} has attained a lot of attention as the limit law is in general different to the case of exogenous sampling. This pointed out that endogeneities can lead to completely new surprising effects and complicate estimators' asymptotic properties. After a concise review of the main findings of \cite{fukasawa}, we shall reveal that similar effects of endogeneity do not arise for multi-scale-type estimators under noise. Thus, the estimation approach is robust against endogeneity of observation times. In light of the core result by \cite{BRHM}, that non-synchronicity effects are asymptotically negligible under microstructure noise in contrast to the non-noisy case, cf.\,\cite{hy3}, this robustness to endogenous sampling in the model with noise does actually not surprise that much. This finding is also in line with recent works by \cite{koike2} and \cite{koike} proving that asymptotic properties of pre-average estimators are not affected by endogenous sampling. In the sequel, we explain the different impact of endogenous sampling in the model with and without noise, respectively.

\subsubsection{Review on Realized Volatility under Endogenous Sampling\label{sec:3.11}}
Assume we have discrete one-dimensional observations $X_{t_i},i\ge 0$, of the process $X$ from \eqref{X} on $[0,T]$ at times $t_i$ given by sequences of stopping times. Consider the martingale approximations $\tilde X_t=X_{t_i}+\int_{t_i}^t \sigma_{t_i}\,dW_s, t\in[t_i,t_{i+1})$ of $X$. The very general and convenient assumption on sampling times by \cite{fukasawa} is of the following nature:
Assume existence of adapted locally bounded left-continuous processes $(a_s)$ and $(b_s)$, such that
for
\begin{align}\label{defendo}\mathcal{X}_{i,n}^k=\inc{\big(\tilde X_{t_{i}}-\tilde X_{t_{i-1}}\big)^k}\end{align}
the following relations with a sequence $n\rightarrow\infty$ are satisfied:
\begin{align}\label{endogcon0}\sum_{t_i\le T}\mathcal{X}_{i,n}^2=\mathcal{O}_{\P}(1)\,,\end{align}
\vspace*{-1cm}

\begin{align}\label{endogcon}\frac{\mathcal{X}_{i,n}^3}{\mathcal{X}_{i,n}^2}=(T/n)^{1/2}b_{t_{i-1}} +\KLEINO_{\P}\big(n^{-1/2}\big),\frac{\mathcal{X}_{i,n}^4}{\mathcal{X}_{i,n}^2}=(T/n)a_{t_{i-1}}^2 +\KLEINO_{\P}\big(n^{-1}\big),\frac{\mathcal{X}_{i,n}^8}{\mathcal{X}_{i,n}^2}=\KLEINO_{\P}\big(n^{-2}\big)\,.\end{align}
In \eqref{defendo} a sequence of sampling schemes is considered such that the right-hand side depends on $n$ which for notational brevity is not highlighted in the times $t_i,i\ge 0$.
Then, realized volatility obeys the stable limit theorem
\begin{align}\label{endogrv}\hspace*{-.175cm}\sqrt{n}\Bigg(\sum_{t_i\le T}\hspace*{-0.1cm}\big(X_{t_{i}}-X_{t_{i-1}}\big)^2\hspace*{-0.05cm}-\hspace*{-0.05cm}\int_0^T\hspace*{-0.1cm}\sigma_s^2\,ds\hspace*{-0.05cm}\Bigg)\hspace*{-0.075cm}\stackrel{st}{\rightarrow}MN\hspace*{-0.075cm}\left(\hspace*{-0.1cm}\sqrt{T}\tfrac23 \int_0^T\hspace*{-0.1cm} b_s\sigma_s\,dW_s,T\hspace*{-0.1cm}\int_0^T \hspace*{-0.1cm}\tfrac23 \big(a_s^2-\tfrac23 b_s^2\big)\sigma_s^2\,ds\hspace*{-0.05cm}\right).\end{align}
For non-random equidistant times $t_i=iT/n,i=0,\ldots,n$, $b_s$ vanishes and $a_s^2=3\sigma_s^2$, such that the random asymptotic variance coincides with $2T\int_0^T\sigma_s^4\,ds$, known from the exogenous sampling case. The most salient point is that the mixed normal limit distribution in \eqref{endogrv} involves a random asymptotic bias. The stable convergence \eqref{endogrv} follows as marginal law from the functional stable convergence
\begin{align}\label{endogrvf}\hspace*{-0.1cm}\sqrt{n}\Bigg(\sum_{t_i\le T}\big(X_{t_{i}}-X_{t_{i-1}}\big)^2\hspace*{-0.05cm}-\hspace*{-0.05cm}\int_0^T\hspace*{-0.1cm}\sigma_s^2\,ds\Bigg)\hspace*{-0.05cm}\stackrel{st}{\rightarrow}\hspace*{-0.05cm}\sqrt{T}\tfrac23 \hspace*{-0.05cm}\int_0^T \hspace*{-0.1cm}b_s\sigma_s\,dW_s\hspace*{-0.05cm}+\hspace*{-0.05cm}\int_0^T\hspace*{-0.1cm}\big(T\tfrac23 \big(a_s^2-\tfrac23 b_s^2\big)\sigma_s^2\big)^{\frac12}\,dB_s\end{align}
with a standard Brownian motion $B$ independent of $W$. As all stable limit theorems in the area of volatility estimation from high-frequency observations of semimartingales, \eqref{endogrvf} is proved based on the theory by \cite{jacod1}.\\
Writing the discretization error of realized volatility
\begin{align*}\sum_{t_i\le T}\big(X_{t_{i}}-X_{t_{i-1}}\big)^2-\int_0^T\sigma_s^2\,ds&=\sum_{t_i\le T}\Big(\big(\tilde X_{t_{i}}-\tilde X_{t_{i-1}}\big)^2-\int_{t_{i-1}}^{t_{i}}\sigma_{s}^2\,ds\Big)+\KLEINO_{\P}(1)\\
&=2\sum_{t_i\le T}\int_{t_{i-1}}^{t_{i}}\big(\tilde X_{s}-\tilde X_{t_{i-1}}\big)d\tilde X_s+\KLEINO_{\P}(1)\,,\end{align*}
five conditions for $\sum_{t_i\le T}\zeta_i^n$, with $\zeta_i^n=2\sqrt{n}\int_{t_{i-1}}^{t_{i}}\big(\tilde X_{s}-\tilde X_{t_{i-1}}\big)d\tilde X_s$, have to be verified according to Theorem 3--1 in \cite{jacod1} in order to conclude \eqref{endogrvf}.\\
The first is satisfied as the $\zeta_i^n$ have vanishing expectation. The following two relations imply the asymptotic law:
\begin{align}&\notag\sum_{t_i\le T}\inc{\big(\zeta_i^n\big)^2}=\sum_{t_i\le T}4 n\inc{\int_{t_{i-1}}^{t_{i}}\big(\tilde X_s-\tilde X_{t_{i-1}}\big)^2d\langle \tilde X\rangle_s}\\
&\quad\notag=\sum_{t_i\le T}n\left(\tfrac23 \inc{\big(\tilde X_{t_{i}}-\tilde X_{t_{i-1}}\big)^4}-\tfrac83\inc{\int_{t_{i-1}}^{t_{i}}\big(\tilde X_s-\tilde X_{t_{i-1}}\big)^3d\tilde X_s}\right)\\
&\quad\notag=\sum_{t_i\le T}T \tfrac23 \inc{\big(\tilde X_{t_{i}}-\tilde X_{t_{i-1}}\big)^2}a_{t_{i-1}}^2+\KLEINO_{\P}(1)\\
&\quad\label{varcond} =\sum_{t_i\le T}T\tfrac23 \inc{\langle \tilde X\rangle_{t_{i}}-\langle \tilde X\rangle_{t_{i-1}}}a_{t_{i-1}}^2+\KLEINO_{\P}(1)=T\int_0^T\tfrac23 \sigma_s^2a_s^2\,ds+\KLEINO_{\P}(1),\end{align}
\begin{align}&\notag \sum_{t_i\le T}\inc{\zeta_i^n \big(\tilde X_{t_{i}}-\tilde X_{t_{i-1}}\big)}=2\sqrt{n}\sum_{t_i\le T}\inc{\int_{t_{i-1}}^{t_{i}}\big(\tilde X_s-\tilde X_{t_{i-1}}\big)d\langle \tilde X\rangle_s}\\
&\quad\notag=2\sqrt{n}\sum_{t_i\le T}\Bigg(\tfrac13 \inc{\big(\tilde X_{t_{i}}-\tilde X_{t_{i-1}}\big)^3}-\inc{\int_{t_{i-1}}^{t_{i}}\big(\tilde X_s-\tilde X_{t_{i-1}}\big)^2d\tilde X_s}\Bigg)\\
&\quad \label{biascond}=\sum_{t_i\le T}\sqrt{T}\tfrac23 \inc{\big(\tilde X_{t_{i}}-\tilde X_{t_{i-1}}\big)^2}b_{t_{i-1}}+\KLEINO_{\P}(1)=\tfrac23 \sqrt{T}\int_0^T b_s\sigma_s^2\,ds+\KLEINO_{\P}(1),\end{align}
together with the two conditions $\sum_{t_i\le T}\inc{\big(\zeta_i^n\big)^4}=\KLEINO_{\P}(1)$ and
$$\sum_{t_i\le T}\inc{\zeta_i^n \big(N_{t_{i+1}}-N_{t_i}\big)}=0$$ for any bounded $\mathcal{F}_t$-martingale $N$ orthogonal to $W$. Contrarily to the exogenous sampling case the local Gaussianity $\sqrt{n}(X_{t_i}-X_{t_{i-1}})\approx N(0,\sigma^2_{t_{i-1}})$ does not apply. We have used It\^{o}'s formula and \eqref{endogcon} to prove the above relations.
For exogenous sampling the term \eqref{biascond} tends to zero, here this non-vanishing covariation is responsible for the asymptotic bias in \eqref{endogrv} and \eqref{endogrvf}.
\subsubsection{Robustness of Multi-Scale Estimators under Endogenous Sampling}
Observe that the crucial difference in \eqref{endogrv} compared to the exogenous case is from the different limits obtained in \eqref{varcond} and \eqref{biascond}. Realized volatility is the sum of squared increments such that third and fourth increments occur in \eqref{varcond} and \eqref{biascond}, respectively. Then, the relations \eqref{endogcon} with processes $b_s,a_s$ determine the asymptotic law.\\
For the multi-scale approach in the model with noise instead, the variance-covariance induced by squared increments is asymptotically negligible and instead cross products of increments over disjoint time segments trigger the (co-)variances of the discretization error, see \eqref{keyendo} in the proofs. Thus similar effects as for \eqref{varcond} and \eqref{biascond} do not occur and generalized It\^{o} isometry implies (co-)variances of the same type as under exogenous observation times.
\begin{assumpsec}\label{disendo}
We have observations at random times $t_i^{(p)},i=0,\ldots,n_p,p=1,\ldots,d$ with $0<\E[n_p]/\E[n_q]<\infty$. We introduce a sequence of sub-filtrations $(\mathcal{F}_t^N)$ of the augmented $(\mathcal{F}_t)$ such that $t_i^{(p)},i=0,\ldots,n_p,p=1,\ldots,d$, are sequences of $(\mathcal{F}_t^N)$-stopping times. For a constant $0<\alpha\le 1/9$, it holds that
\begin{align}\label{eqdis4}
\delta_n&=\sup_{(i,p)}{\left(\left(t_i^{(p)}-t_{i-1}^{(p)}\right),t_0^{(p)},T-t_{n_p}^{(p)}\right)}~\;\,=\mathcal{O}_{\P}\left(\sup_p(n_p)^{-\frac{8}{9}-\alpha}\right)~.\end{align}
Assume stochastic convergence of the sequences in Assumption \ref{aqvt2} (i). When the indicator functions in Assumption \ref{aqvt2} are replaced by $\P\big({t_p^+(T_j)=t_{p'}^+(T_j)}|\mathcal{F}_{T_j}^N\big)$ and analogously for the other sets, assume convergence of the respective series. We use the same notation for the limit objects as above.
\end{assumpsec}
\begin{cor}\label{endocor}On the Assumptions \ref{eff} and \ref{disendo}, the generalized multi-scale estimator \eqref{gms} with $M_{N}=c\sqrt{N}$ and weights \eqref{weights} obeys the multivariate stable central limit theorem \eqref{gmsclt} with asymptotic variance-covariance matrix \eqref{covgms}.
\end{cor}
\subsection{The Feasible Multivariate Limit Theorem\label{sec:3.3}}
A remaining step towards a feasible asymptotic distribution theory allowing for confidence is to provide a consistent estimator for the asymptotic variance-covariance matrix \eqref{covgms}.
In the vein of \cite{bibinger2}, we construct a consistent estimator in the general non-synchronous framework, following a histogram-type approach.
\begin{prop}\label{covest}
On the assumptions imposed in Theorem \ref{procovgms}, the estimator
\begin{align}\label{covestgms}\widehat{\Acov}&=4\,c\sum_{j=1}^{K_N}\Big(\hat\Sigma_{\frac{(j-1)T}{K_N}}^D\otimes \hat\Sigma_{\frac{(j-1)T}{K_N}}^D\Big)\mathcal{Z}\frac{D^{\alpha}_N(T)}{K_N}+2\,c^{-3}\widehat{\big(\mathbf{H}\otimes\mathbf{H}\big)}^*\mathcal{Z}\\
&\notag\quad +c^{-1}\mathfrak{M}^{\alpha}\sum_{j=1}^{K_N}\mathcal{Z}\Bigg(\Bigg(\frac{(S^N(T))^{(pq)}}{K_N}\hat{\mathbf{H}}^{(pq)}\Bigg)_{1\le p,q\le d}\otimes \hat\Sigma_{\frac{(j-1)T}{K_N}}^S\Bigg)\mathcal{Z}+c^{-1}\mathfrak{N}_2^{\alpha}\widehat{\big(\tilde{\mathbf{H}}^{\otimes}\big)}\mathcal{Z}\,,\end{align}
gives a consistent estimator of \eqref{covgms}. The estimators $\hat\Sigma_{\frac{(j-1)T}{K_N}}^D$ and $\hat\Sigma_{\frac{(j-1)T}{K_N}}^S$ are multi-scale estimators over bins $[D_{j-1}^N,D_j^N]$ and $[(S_{j-1}^N)^{(pq)},(S_j^N)^{(pq)}],j=1,\ldots,K_N$, rescaled with $\Delta_j D^N$ and $\big(\d_j S^N\big)^{(pq)}$, respectively, with multi-scale frequencies $M_N(j),1\le j\le K_N$. Here we use
\begin{align*}D_j^N=\inf{\left\{t\in[0,T]|D_N^{\alpha}(t)\ge jD_N^{\alpha}(T)/K_N\right\}},0\le j\le K_N\,,~\Delta_j D^N=D_j^N-D_{j-1}^N\,,\end{align*}
\begin{align*}(S_j^N)^{(pq)}=\inf{\left\{t\in[0,T]|(S^N(t))^{(pq)}\ge j(S^N(T))^{(pq)}/K_N\right\}},0\le j\le K_N,\end{align*}
\begin{align*}\big(\Delta_j S^N\big)^{(pq)}=(S_j^N)^{(pq)}-(S_{j-1}^N)^{(pq)},1\le j\le K_N,\end{align*}
with $D_N^{\alpha}$ the right-hand side of \eqref{wlasa} and $S^N$ from \eqref{dr} which hinge on the weights and the observation times. In the ex-post estimation we can use the available observation times. With estimators of the noise variance-covariance, e.g.\,
\begin{align*}\hat{\mathbf{H}}^{(pp)}=-n_p^{-1}\sum_{i=1}^{n_p-1}\d_i Y^{(p)}\d_{i+1}Y^{(p)}\,,\end{align*}
and for synchronous observations analogous for non-diagonal entries, we readily obtain estimates for the remaining terms above.
The number of bins $K_N$ is chosen such that $K_NN^{-\nicefrac{1}{3}}\rightarrow 0, K_N\rightarrow\infty$ as $N\rightarrow\infty$. A possible choice is $K_N=cN^{1/5}$ and $M_N(j)=c^{5/4}N^{3/5}$. We derive the feasible multivariate central limit theorem
\begin{align}N^{1/4}\widetilde{\Acov}^{-1/2}\Big(\widehat{[X,X]}_T^{(multi)}-\int_0^T\Sigma_s\,ds\Big)\stackrel{st}{\rightarrow} N(0, \mathcal{Z})\,,\end{align}
with $\widehat{\Acov}=\widetilde{\Acov}\mathcal{Z}$.\end{prop}
\begin{remark}
The feasible limit theorem remains valid when relaxing Assumptions of Theorem \ref{procovgms} on existence of the limit $\mathfrak{G}$, since every subsequence of \eqref{wlasa} has an in probability converging subsequence, see the discussion at the end of p.\,1411 in \cite{zhangmykland} for analogous reasoning and more details.
\end{remark}
The estimator \eqref{covestgms} simplifies in many cases, i.\,e.\,the completely non-synchronous setup according to Corollary \ref{corcovgms}. An estimator for \eqref{covms} in the synchronous case is contained in \eqref{covestgms}. It is natural for multivariate feasible limit theorems that the known non-invertible matrix $\mathcal{Z}$ remains in the limit law.
\section{An Application to Conditional Independence Testing\label{sec:4}}
This section is devoted to the design of a statistical test in order to investigate if the correlation of two assets is only induced by a factor to which both are correlated. For portfolio modeling and management, information about such relations can provide valuable information and access to a new angle on the covariance structure. Conclusions that significant integrated covolatilities between high-frequency assets are fully explained by their dependence on a joint factor or another asset, respectively, facilitate dimension reduction of covariance matrix estimation which is particularly important when considering multivariate limit theorems with variance-covariance matrices of dimension \eqref{dimension}. Consequently, each relation \eqref{teststatistic} equal to zero reduces the required estimates in \eqref{covestgms}. Beyond this practical implication, a relation \eqref{teststatistic} equal to zero reveals knowledge about the dependence structure useful e.g.\,for default contagion as well as for many other economic applications. For instance, we can think of two observed asset processes $X_1$ and $X_2$ listed within one index $Z$ being conditionally on $Z$ independent. To put it the other way round, pairs which are not conditionally independent exhibit significant covariance that carries information about the direct mutual influence. We understand independence here in terms of orthogonal quadratic covariation processes and test for zero integrated covolatility -- so the term `independence' is used here for a simple illustrative phrasing. $X_1$ and $X_2$ are orthogonally decomposed in the sum of $Z$ and a process independent of $Z$. The constants $\rho^{X_1},\rho^{X_2}$ quantify the degree of dependence on $Z$.
\begin{align} X_1=\rho^{X_1}\,Z+Z^{\bot}\,,\,X_2=\rho^{X_2}\,Z+Z^{\dagger}~~\mbox{with}~[Z,Z^{\bot}]\equiv 0\,,\,[Z,Z^{\dagger}]\equiv 0\,.\end{align}
With $[X_1,X_2]\equiv 0$ for two semimartingales $X_1,X_2$ we express that $[X_1,X_2]_s=0$ for all $s\in[0,T]$. For the conditional independence hypothesis, we set \begin{align}\mathds{H}_0: [Z^{\bot},Z^{\dagger}]_T=0\,.\end{align}
Essentially, we do not distinguish between pairs for which the orthogonal parts are uncorrelated on the whole line and pairs for which this correlation process integrates to zero. Our focus is on a resulting zero quadratic covariation over $[0,T]$.\\
A suitable test statistic to decide whether we reject $\mathds{H}_0$ or not is
\begin{align}\label{teststatistic}\mathfrak{T}(X_1,X_2,Z)=[X_1,Z]_T[X_2,Z]_T-[X_1,X_2]_T[Z,Z]_T\,,\end{align}
which is zero under $\mathds{H}_0$.\\ In our high-frequency framework we can estimate the single integrated (co-)volatilities via the approaches considered in the preceding sections. The vital point is to deduce the asymptotic distribution of the estimated version
\begin{align}\hat {\mathfrak{T}}_n=\widehat{[X_1,Z]}_T^{(multi)}\widehat{[X_2,Z]}_T^{(multi)}-\widehat{[X_1,X_2]}^{(multi)}_T\widehat{[Z,Z]}^{(multi)}_T\,,\end{align}
where $\widehat{[\,\cdot\,]}_T^{(multi)}$ stands for one of the aforementioned estimators \eqref{ms} or \eqref{gms}. This test statistic, though based on the simple function $g(x,y,u,v)=xy-uv$, is more complex to analyze than linear combinations, since we face products of our estimators. Therefore, the asymptotic law of \eqref{teststatistic} is not directly obtained from Theorem \ref{procovms} or Theorem \ref{procovgms}, respectively. In lieu of determining the distribution of the test statistic, we apply the $\Delta$-method for stable convergence. Here, the methodology is similar to the prominent propagation of error concept from experimental science based on a simple Taylor expansion. For each quadratic covariation, the estimation error gets small for large $n$ and hence we can profit when we taylor the underlying function $g$. Indeed, this will give us the leading term of the variance of $\hat{\mathfrak{T}}_n$:
\begin{align}&\mathfrak{T}-\hat{ \mathfrak{T}}_n=[X_2,Z]_T\left([X_1,Z]_T-\widehat{[X_1,Z]}_T^{(multi)}\right)+[X_1,Z]_T\left([X_2,Z]_T-\widehat{[X_2,Z]}_T^{(multi)}\right)\\
&\notag -[X_1,X_2]_T\left([Z,Z]_T-\widehat{[Z,Z]}_T^{(multi)}\right)-[Z,Z]_T\left([X_1,X_2]_T-\widehat{[X_1,X_2]}_T^{(multi)}\right)+\mathcal{O}_{\P}\big(n^{-\frac12}\big)\,.\end{align}
The asymptotic variance of the test statistic is random as a linear combination of the unknown quadratic (co-)variations and entries of the asymptotic variance-covariance matrix. We have consistently estimated the latter in Proposition \ref{covest}. Denote by $\AVAR(U)$ and $\Acov(U,V)$ in the sequel asymptotic variances and covariances of one-dimensional random variables $U,V$. An elementary calculation yields
\begin{align*}
\AVAR(\hat {\mathfrak{T}}_n)&=[X_2,Z]_T^2\AVAR\big(\widehat{[X_1,Z]}_T^{(multi)}\big)+[X_1,Z]_T^2\AVAR\big(\widehat{[X_2,Z]}_T^{(multi)}\big)\\
&\quad+[X_1,X_2]_T^2\AVAR\big(\widehat{[Z,Z]}_T^{(multi)}\big)+[Z,Z]_T^2\AVAR\big(\widehat{[X_1,X_2]}_T^{(multi)}\big)\\
&\quad+2\,[Z,Z]_T[X_1,X_2]_T\Acov\big(\widehat{[X_1,X_2]}_T^{(multi)},\widehat{[Z,Z]}_T^{(multi)}\big)\\
&\quad +2\,[X_1,Z]_T[X_2,Z]_T\Acov\big(\widehat{[X_1,Z]}_T^{(multi)},\widehat{[X_2,Z]}_T^{(multi)}\big)\\
&\quad-2\,[X_1,Z]_T[Z,Z]_T\Acov\big(\widehat{[X_1,X_2]}_T^{(multi)},\widehat{[X_2,Z]}_T^{(multi)}\big)\\
&\quad -2\,[X_2,Z]_T[Z,Z]_T\Acov\big(\widehat{[X_1,X_2]}_T^{(multi)},\widehat{[X_1,Z]}_T^{(multi)}\big)\\
&\quad-2\,[X_1,X_2]_T[X_1,Z]_T\Acov\big(\widehat{[X_2,Z]}_T^{(multi)},\widehat{[Z,Z]}_T^{(multi)}\big)\\
&\quad -2\,[X_1,X_2]_T[X_2,Z]_T\Acov\big(\widehat{[X_1,Z]}_T^{(multi)},\widehat{[Z,Z]}_T^{(multi)}\big).
\end{align*}
Inserting consistent estimators for the asymptotic (co-)variances above, we obtain with our multivariate stable central limit theorem that
\begin{align}n^{\frac14}\Big(\widehat{\AVAR(\hat {\mathfrak{T}}_n)}\Big)^{-1/2}\hat {\mathfrak{T}}_n\stackrel{st}{\rightarrow}N(0,1)\,,\end{align}
or with scaling $N^{1/4}$ for non-synchronous observations, under $\mathds{H}_0$ what gives an asymptotic distribution free test.\\
The role of $Z$ in the model can be also some macro variable that is either known or can be estimated with faster rate of convergence which simplifies the terms above. For regularly observed high-frequency data without noise the same kind of test can be constructed using the realized volatility matrix and the faster rate $n^{1/2}$ is attained.

\section{An Empirical Example\label{sec:5}}
We survey our methods in an application study on NASDAQ intra-day trading data, reconstructed from first-level order book data, from August 2010. We consider a sample portfolio with $5$ assets, namely Apple (AAPL), Microsoft (MSFT), Oracle (ORCL), Exxon Mobil Corporation (XOM) and Pfizer (PFE). Traded prices are recorded at non-synchronous times and market microstructure noise is clearly indicated such that we suppose the model from Assumption \ref{dis4}. We quantify the integrated volatility matrix over the whole month (where we discard over-night returns) and for the first trading day, 2010/08/02, respectively, using generalized multi-scale estimates \eqref{gms} with weights \eqref{cubweights} and pairwise refresh times. The complete variance-covariance matrix of the estimates is quantified. 
For a $d$-dimensional portfolio, the number of free entries of this symmetric variance-covariance matrix is given by
\begin{align}\label{dimension}\frac{1}{2}\frac{d(d+1)}{2}\left(\frac{d(d+1)}{2}+1\right)=d+3\binom{d}{4}+3\cdot 2\binom{d}{3}+4\binom{d}{2}\,.\end{align}
The left-hand illustration is derived as $d(d+1)/2$ different entries of the symmetric integrated volatility matrix are estimated which leads to a $(d(d+1)/2)\times(d(d+1)/2) $-dimensional variance-covariance matrix of estimates which is symmetric again. 
In Table \ref{tabdat1}, we list the estimates for the integrated volatility matrices $\pm$ estimated standard deviations. The variance-covariance matrices of these estimates based on \eqref{covestgms} and the numbers of pairwise refresh times are listed in Table \ref{tabdat2}. The bottom line is that involving covariances of estimates is indispensable when facing questions for multivariate portfolio management. The estimated quadratic variation of a sum of all five assets is $(41.57\pm0.26)\cdot 10^{-3}$ for 2010/08 and $(129.22\pm6.88)\cdot 10^{-5}$ for 2010/08/02. The risk of the estimated volatilities for these portfolios, $6.93\cdot 10^{-8}$ and $47.42\cdot 10^{-10}$, is mainly induced by covariances (6.34/42.91), whereas the trace of the variance-covariance matrix, i.\,e.\,the sum of estimated variances, is much smaller. If one would mistakenly act as if the estimators were uncorrelated, this leads to a tremendous underestimate of uncertainty.\\
We perform the test from Section \ref{sec:4} to investigate three hypotheses: if MSFT and ORCL have a zero covariation conditional on PFE; ORCL and PFE conditional on MSFT and MSFT and PFE conditional on ORCL. 
We obtain the following $p$-values as test results  
\begin{align*}
&~p=0.00;0.014;0.064~~\text{(2010/08)}\,,\\
&~p=0.00;~0.96~;~0.21~~~\text{(2010/08/02)}.\end{align*}
Further tests if MSFT and ORCL have zero covariation conditional on the sum of all assets yield $p$-values 0.13 for 2010/08 and 0.23 for 2010/08/02. Tests for ORCL and PFE conditional on the sum of the three other assets yield $p$-values 0.19 and 0.99.\\
In conclusion, this empirical evidence suggests that MSFT and ORCL have some dependence not explained by a common macro factor influencing all NASDAQ assets. On the contrary, we can not reject this for several other combinations. Some differences between 2010/08 and 2010/08/02 give an heuristic that the portfolio dependence structure is not completely persistent. Though there are some limitations where the additive noise model does not perfectly fit the stylized facts of the considered high-frequency data as discreteness of returns and zero returns, the approaches developed in this research area and advancements of this article provide reliable tools to quantify risk measures from high-frequency asset prices and determine confidence intervals for the estimates.

\begin{table}[t]
\begin{center}
\begin{tabular}{|l||c|c|c|c|c|}
\hline
 $\widehat{[ X,X]}_T^{(multi)}$~ & AAPL & MSFT & ORCL & XOM & PFE \\
 \hline\hline
 AAPL & $2.82\pm 0.01$& $1.23\pm 0.02$&$1.66\pm 0.03$&$0.63\pm 0.01$&$0.99\pm 0.03$\\
 MSFT & \cellcolor{gray} &$2.87\pm 0.03$& $2.05\pm 0.05$&$0.94\pm 0.02$&$1.25\pm 0.07$\\
 ORCL & \cellcolor{gray} & \cellcolor{gray} &$3.81\pm 0.06$&$1.31\pm 0.03$&$1.63\pm 0.10$\\
 XOM & \cellcolor{gray} & \cellcolor{gray} & \cellcolor{gray} &$2.36\pm 0.01$&$1.25\pm 0.04$\\
 PFE  & \cellcolor{gray} & \cellcolor{gray} & \cellcolor{gray}& \cellcolor{gray}& $3.83\pm 0.09$\\
 \hline
\end{tabular}\\
\begin{tabular}{|l||c|c|c|c|c|}
\hline
 $\widehat{[ X,X]}_T^{(multi)}$~ & AAPL & MSFT & ORCL & XOM & PFE \\
 \hline\hline
 AAPL & $8.16\pm 0.10$& $3.42\pm 0.20$&$3.74\pm 0.20$&$2.28\pm 0.11$&$2.47\pm 0.26$\\
 MSFT & \cellcolor{gray} &$12.43\pm 0.38$& $7.43\pm 0.37$&$2.95\pm 0.20$&$4.69\pm 0.47$\\
 ORCL & \cellcolor{gray} & \cellcolor{gray} &$11.90\pm 0.39$&$3.73\pm 0.22$&$2.90\pm 0.49$\\
 XOM & \cellcolor{gray} & \cellcolor{gray} & \cellcolor{gray} &$6.42\pm 0.12$&$1.75\pm 0.29$\\
 PFE  & \cellcolor{gray} & \cellcolor{gray} & \cellcolor{gray}& \cellcolor{gray}& $19.59\pm 0.70$\\
 \hline
\end{tabular}
\end{center}
\caption{\label{tabdat1}Estimates for the integrated volatility matrix ($\cdot$ $10^3$) 2010/08 (top) and ($\cdot$ $10^5$) 2010/08/02 (bottom).}
\end{table}
\newpage
\begin{landscape}
\small
\enlargethispage*{5cm}
\vspace*{-1.5cm}

\begin{longtable}{l|ccccccccccccccc}
&[A,A]&[A,M]&[A,O]&[A,X]&[A,P]&[M,M]&[M,O]&[M,X]&[M,P]&[O,O]&[O,X]&[O,P]&[X,X]&[X,P]&[P,P]\\
\hline
[A,A]& 4.01& 4.53& 5.47& 1.55& 5.26& 6.97& 6.75& 4.58& 6.92& 9.33& 5.15& 8.28& 3.02& 3.46& 10.89\\\relax
[A,M]& & 4.91& 4.63& 1.89& 2.72& 7.48& 7.74& 3.33& 2.72& 8.74& 2.88& 9.30& 1.22&	2.88& 12.23\\\relax
[A,O]&&& 7.21&	1.85&	4.50	&5.85	&2.89&	2.91&	7.26&	11.62&4.70&	9.61&	2.03	&2.80	&3.85\\\relax
[A,X]&&&& 2.83	&2.85	&6.15	&5.98	&3.44	&6.71	&8.03	&1.94	&7.49	&3.97	&4.59	&2.96\\\relax
[A,P]&&&&& 7.03	&4.30	&5.63	&2.72	&8.96&	6.71&	0.63&	11.34&	1.37&	4.71&	13.72\\\relax
[M,M]&&&&&&	13.31	&12.29	&5.81&	8.18&	8.06&	3.33	&6.72&	1.21&	2.90&	5.33\\\relax
[M,O]&&&&&&& 12.36	&5.86	&10.18	&15.29	&5.71	&11.07	&1.35	&3.61	&6.32\\\relax
[M,X]&&&&&&&& 5.24	&6.95&	6.03&	5.48&	6.15&	2.72&	4.51&	6.54\\\relax
[M,P]&&&&&&&&& 13.86	&6.82	&4.02	&14.42	&1.82&	6.25&	12.47\\\relax
[O,O]&&&&&&&&&&	21.04&	7.85&	12.93&	1.46&3.98&	7.50\\\relax
[O,X]&&&&&&&&&&&	6.67&	9.08&	2.97&	4.74&	7.98\\\relax
[O,P]&&&&&&&&&&&&17.28&	1.72&	7.63&	15.67\\\relax
[X,X]&&&&&&&&&&&&&4.90&	4.92&	13.11\\\relax
[X,P]&&&&&&&&&&&&&&8.20&	12.23\\\relax
[P,P]&&&&&&&&&&&&&&&29.98\\
\hline
\addtocounter{table}{-1}
\end{longtable}
\begin{longtable}{l|ccccccccccccccc}
&[A,A]&[A,M]&[A,O]&[A,X]&[A,P]&[M,M]&[M,O]&[M,X]&[M,P]&[O,O]&[O,X]&[O,P]&[X,X]&[X,P]&[P,P]\\
\hline
[A,A]& 1.55& 2.01& 1.75& 0.77& 0.65& 1.52& 1.50& 0.41& 0.62& 1.91& 0.73& 1.01& 0.27& 0.31& 2.05\\\relax
[A,M]& & 2.64& 2.32& 0.92& 1.78& 5.67& 4.32& 2.46& 4.27& 3.12& 1.79& 2.89& 0.56&	1.08& 3.05\\\relax
[A,O]&&& 2.35&	1.12&	1.22	&2.96	&4.05&	1.42&	2.35&	5.21&	2.47&	3.54&	0.70	&0.77	&2.63\\\relax
[A,X]&&&& 0.86	&0.68	&1.84	&1.85	&1.38	&0.95	&2.18	&1.53	&1.05	&1.04	&0.90	&2.14\\\relax
[A,P]&&&&& 3.29	&2.52	&1.80	&0.91	&5.57&	1.81&	0.61&	5.64&	0.28&	3.29&	7.99\\\relax
[M,M]&&&&&&	13.78	&9.94	&4.79&	11.05&	5.72&	2.73	&4.59&	0.68&	1.74&	4.66\\\relax
[M,O]&&&&&&& 7.27	&4.1	&5.86	&9.65	&4.33	&6.79	&0.84	&1.64	&3.03\\\relax
[M,X]&&&&&&&& 2.3	&2.68&	3.15&	2.46&	2.05&	1.35&	2.35&	2.53\\\relax
[M,P]&&&&&&&&& 9.21	&2.60	&0.90	&9.49	&0.34&	4.56&	15.32\\\relax
[O,O]&&&&&&&&&&	12.21&	5.63&	5.51&	0.92&	1.14&	2.84\\\relax
[O,X]&&&&&&&&&&&	2.63&	1.78&	1.72&	3.17&	1.43\\\relax
[O,P]&&&&&&&&&&&&9.88&	0.33&	5.18&	7.03\\\relax
[X,X]&&&&&&&&&&&&&1.62&	1.08&	2.2\\\relax
[X,P]&&&&&&&&&&&&&&3.84&	5.42\\\relax
[P,P]&&&&&&&&&&&&&&&34.16\\
\hline
\caption{\label{tabdat2}Estimated asymptotic covariance matrix ($\cdot$ $10^{8}$), 2010/08 (top) and ($\cdot$ $10^{10}$) 2010/08/02 (below), A=AAPL, M=MSFT, O=ORCL, X=XOM, P=PFE.}
\end{longtable}

\end{landscape}
\newpage\normalsize\noindent
\subsection*{Acknowledgements}
The authors would like to thank
Dacheng Xiu, Ruey Tsay and Dan Christina Wang for inspiring discussions. We also thank Johannes Schmidt-Hieber and Till Sabel for a discussion on the relation of estimators via quadratic forms for noise smoothing who have worked on the transformations between several of these estimators.
\\
Markus Bibinger was supported by a fellowship within the Postdoc-Programme of the German Academic Exchange Service (DAAD) and the CRC 649 `Economic Risk' at Berlin, supported by the Deutsche Forschungsgemeinschaft (DFG), and also from the Stevanovich Center for
Financial Mathematics at the University of Chicago.
Per Mykland was supported under National Science Foundation grant SES 11-24526. We gratefully acknowledge this financial support.
\\
The NASDAQ trading data used in Section \ref{sec:5} has been reconstructed from the order book with LOBSTER provided by the high-frequency econometrics team at the Chair of Econometrics, School of Business and Economics, Humboldt-Universität zu Berlin, Germany.

\appendix
\setcounter{equation}{0}
\renewcommand{\theequation}{\thesection.\arabic{equation}}%
\vspace*{1cm}
\begin{center}
\Large%

\textbf{APPENDIX:\ PROOFS}

\renewcommand{\baselinestretch}{1.0}
\normalsize%

\end{center}


\section{Preliminaries\label{sec:9.1}}
The local boundedness condition in Assumption \ref{eff} can be strengthened to uniform boundedness on $[0,T]$ by a localization procedure carried out in \cite{jacodlecture},  Lemma 6.\,6 of Section 6\,.3.
Let $C$ be a generic constant and denote $\d_i W=W_{t_i}-W_{t_{i-1}},i=1,\ldots,n$, for the Brownian motion $W$ driving the SDE with solution $X$ in \eqref{X} and $\d_i \sigma=\sigma_{t_i}-\sigma_{t_{i-1}}$. Consider some norm $\|\,\cdot\,\|$, e.\,g.\,the euclidean norm, on $\mathds{R}^d$. Suppose Assumption \ref{eff} holds. By several applications of the Burkholder-Davis-Gundy and Hölder inequality one can obtain the following estimates:
\begin{subequations}
\begin{align}\label{prel1}\E\left[\|\d_i X\|^2+\|\d_i W\|^2\big|\mathcal{F}_{t_{i-1}}\right]\le C n^{-1},~~~\E\left[\|\d_i \sigma\|^2\big|\mathcal{F}_{t_{i-1}}\right]\le C n^{-1}\,,\\ \label{prel2}\E\left[\|\d_i X-\sigma_{t_{i-1}}\d_i W\|^2\big|\mathcal{F}_{t_{i-1}}\right]\le C n^{-2}\,,\end{align}\end{subequations}
for equidistant observation schemes $t_i=iT/n$. For general synchronous sampling \eqref{prel1} and \eqref{prel2} remain valid when replacing $n$ by $\delta_n^{-1}$ with $\delta_n=\sup_i{\left(t_i-t_{i-1}\right)}$. The estimates \eqref{prel1} and \eqref{prel2} are proven in \cite{jacodlecture}, among others. They are used repeatedly in the analysis below.
We write $a_n\asymp^p b_n$ if $a_n=\mathcal{O}_{\P}(b_n)$ and $b_n=\mathcal{O}_{\P}(a_n)$ and express analogously $a_n\asymp b_n$ for $a_n=\mathcal{O}(b_n)$ and $b_n=\mathcal{O}(a_n)$.\\
A summary including the elements of matrix algebra which are heavily used throughout the proofs can be found in Sections 10.1, 10.2 and 11.2 of \cite{magnus}.
Let us calculate next the asymptotic variance-covariance matrix of the realized volatility matrix in \eqref{cltrv} which serves as well as preparation for the proofs below. Denote by $Z_i\in\R^d,i=0,\ldots,n$, independent standard normally distributed random vectors. We apply the rule $\vec(ABC)=(C^\top \otimes A)\vec(B)$ for matrices $A,B,C$ frequently below. The multivariate stable central limit theorems are proved based on Theorem 3--1 of \cite{jacod1}. The limiting variance-covariance matrix in \eqref{cltrv} is random and obtained, as discussed in Section \ref{sec:3.11}, as the stochastic limit of the sum of conditional variance-covariance matrices.  We find that
\begin{align*}
&\sum_{i=1}^n\cov\Big(\vec\Big(\d_i X(\d_i X)^{\top}\Big)\Big|\mathcal{F}_{\frac{(i-1)T}{n}}\Big)\\
&\hspace*{1cm}\asymp^p\sum_{i=1}^n\cov\Big(\vec\Big(\sqrt{\frac{T}{n}}\Sigma_{\frac{(i-1)T}{n}}^{1/2}\big(Z_i(Z_i)^{\top}\big)\sqrt{\frac{T}{n}}\Sigma_{\frac{(i-1)T}{n}}^{1/2}\Big)\Big|\mathcal{F}_{\frac{(i-1)T}{n}}\Big)\\
&\hspace*{1cm}= \frac{T^2}{n^2}\sum_{i=1}^n\cov\Bigg(\Big(\Sigma_{\frac{(i-1)T}{n}}^{1/2}\otimes \Sigma_{\frac{(i-1)T}{n}}^{1/2}\Big)\vec\big(Z_iZ_i^{\top}\big)\Big|\mathcal{F}_{\frac{(i-1)T}{n}}\Bigg)\\
&\hspace*{1cm}=\frac{T^2}{n^2}\sum_{i=1}^n \Big(\Sigma_{\frac{(i-1)T}{n}}^{1/2}\otimes \Sigma_{\frac{(i-1)T}{n}}^{1/2}\Big)\mathcal{Z}  \Big(\Sigma_{\frac{(i-1)T}{n}}^{1/2}\otimes \Sigma_{\frac{(i-1)T}{n}}^{1/2}\Big)\\
&\hspace*{1cm}=\frac{T}{n}\sum_{i=1}^n \Big(\Sigma_{\frac{(i-1)T}{n}}\otimes \Sigma_{\frac{(i-1)T}{n}}\Big)\frac{T}{n}\mathcal{Z}\asymp^p \frac{T}{n}\int_0^T \Big(\Sigma_{s}\otimes \Sigma_{s}\Big)\mathcal{Z}\,ds\,.
\end{align*}
All other ingredients required to conclude \eqref{cltrv} by Theorem 3--1 of \cite{jacod1}, cf.\,Section \ref{sec:3.11} above, are easily established here and we skip the details.\\
We use analogous transforms for computing terms of the form
 \[\cov\big(\vec\big(\d_j X\otimes (\d_l X)^{\top}\big)\big)=\cov\big(\vec\big(\d_j X(\d_l X)^{\top}\big)\big)\]
frequently below without repeating each step.
\section{Proofs of Section \ref{sec:2}}
\subsection{Proof of Theorem \ref{thm:comparison-raw}}
For the proof that
\begin{align*}\widehat{\left[ X, X\right]}_T^{(multi)}-\widehat{\left[ X, X\right]}_T^{(kernel)}=-4\mathbf{H}+\KLEINO_{\P}\left(n^{-\frac{1}{4}}\right)\end{align*}
if $\mathfrak{K}^{\prime\prime}=h$ in \eqref{weights}, it suffices to focus on the first-order term of the weights. Transforming \eqref{ms} yields
\begin{align*}\sum_{i=1}^{M_n}\frac{\alpha_{i}}{i}\sum_{j=i}^n\d_j^iY(\d_j^iY)^{\top}&=\sum_{i=1}^{M_n}\alpha_i\left(\sum_{j=2}^n\sum_{l=1}^{i \wedge (j-1)}\left(1-\frac{l}{i}\right)\big(\d_jY(\d_{j-l}Y)^{\top}+\d_{j-l}Y(\d_jY)^{\top}\big)\right)\\
&\hspace*{1cm}+\sum_{j=1}^n\d_j Y\d_j Y^{\top}-R_n\\
&=\sum_{l=1}^{M_n}\sum_{j=l+1}^n\sum_{i=l}^{M_n}\alpha_i\left(1-\frac{l}{i}\right)\big(\d_jY(\d_{j-l}Y)^{\top}+\d_{j-l}Y(\d_jY)^{\top}\big)\\
&\hspace*{1cm}+\sum_{j=1}^n\d_j Y(\d_j Y)^{\top}-R_n\,.
\end{align*}
The term $R_n$ induced by end-effects
\begin{align*}\hspace*{-.05cm}&\sum_{i=1}^{M_n}\alpha_i\left(\sum_{j=1}^{i-1}\left(\frac{i-j}{i}\d_j Y(\d_j Y)^{\top}+\sum_{l=1}^{(j-1)\wedge 1}\frac{i-j}{i}\big(\d_jY(\d_{j-l}Y)^{\top}+\d_{j-l}Y(\d_jY)^{\top}\big)\right)\right.\\
\hspace*{-.05cm}&\left.+\hspace*{-.1cm}\sum_{j=n-i+2}^{n}\hspace*{-.1cm}\left(\frac{i-n+j-1}{i}\left(\d_j Y(\d_j Y)^{\top}\hspace*{-.1cm}+\hspace*{-.1cm}\sum_{l=1}^{i\wedge (n-j)}\hspace*{-.1cm}\big(\d_jY(\d_{j-l}Y)^{\top}+\d_{j-l}Y(\d_jY)^{\top}\big)\right)\hspace*{-.1cm}\right)\hspace*{-.1cm}\right)\end{align*}
has an expectation by noise:
\begin{align*}2\,{\bf{H}}\,\sum_{i=1}^{M_n}\alpha_i\left(\sum_{j=1}^{i-1}\frac{i-j}{i}-\sum_{j=2}^{i-1}\frac{i-j}{i}+\sum_{j=n-i+1}^{n-1}\frac{i-n+j}{i}-\sum_{j=n-i+1}^{n-2}\frac{i-n+j}{i}\right)=4\,{\bf{H}}\,.\end{align*}
The variance-covariance matrix of this term is asymptotically negligible what can be shown with standard bounds. For the main term above, we can detach the inner sum and find that
\begin{align*}\sum_{i=l}^{M_n}\alpha_i\left(1-\frac{l}{i}\right)&=\sum_{i=l}^{M_n}\frac{i}{M_n^2}\frac{(i-l)}{i}\mathfrak{K}^{\prime\prime}\left(\frac{i}{M_n}\right)+\KLEINO\left(n^{-\frac{1}{4}}\right)\\
&=\int_{l/M_n}^1\mathfrak{K}^{\prime\prime}(x)\left(x-\frac{l}{M_n}\right)\,dx+\KLEINO\left(n^{-\frac{1}{4}}\right)=\mathfrak{K}\left(\frac{l}{M_n}\right)+\KLEINO\left(n^{-\frac{1}{4}}\right)\,,\end{align*}
by partial integration under the restrictions made on $\mathfrak{K}$. This yields the form \eqref{kernel} of the transformed kernel estimator and our claim. That the integral approximation does not harm the above equality up to the $\KLEINO\left(n^{-1/4}\right)$-term, can be seen by the estimate
\begin{align*}\int_{i/M_n}^{(i+1)/M_n}\left|f(x)-f\left(\frac{i}{M_n}\right)\right|\,dx\le \int_{i/M_n}^{(i+1)/M_n}C\,\left|x-\frac{i}{M_n}\right|\,dx\le C\, M_n^{-2}\end{align*}
with generic constant $C$, $i\ge l$, for the Lipschitz function
\begin{align*}f(x)=\mathfrak{K}^{\prime\prime}(x)\left(x-\frac{l}{M_n}\right)\end{align*}
on the compact support $[0,1]$, where Lipschitz continuity is ensured by the preconditioned continuous differentiability.\\
The extension of the equivalence from $H_n=M_n$ in Theorem \ref{thm:comparison-raw} to asymptotically of the same (optimal) order
follows directly, by inserting the minimum in the transformations above, and by elementary bounds for the remainder.$\hfill\Box$

\subsection{Proof of Theorem \ref{procovms}}
Decompose the multi-scale estimator \eqref{ms}
\begin{align*}\sum_{i=1}^{M_n}\frac{\alpha_{i}}{i}\sum_{j=i}^{n}\d^{i}_jY(\d^{i}_j Y)^{\top}&=\sum_{i=1}^{M_n}\frac{\alpha_{i}}{i}\sum_{j=i}^{n}\d^{i}_j X(\d^{i}_j  X)^{\top}+\sum_{{i}=1}^{M_n}\frac{\alpha_{i}}{i}\sum_{j=i}^{n}\d^{i}_j\epsilon(\d^{i}_j \epsilon)^{\top}\\
&+\sum_{i=1}^{M_n}\frac{\alpha_{i}}{i}\sum_{j=i}^{n}\d^{i}_j X(\d^{i}_j \epsilon)^{\top}+\sum_{i=1}^{M_n}\frac{\alpha_{i}}{i}\sum_{j=i}^{n}\d^{i}_j\epsilon(\d^{i}_j  X)^{\top}\,,\end{align*}
with $\d_i^j\epsilon=\epsilon_j-\epsilon_{j-i}$, in a signal part, a noise part and cross terms which are uncorrelated. We analyze the variance-covariance matrices of the signal, noise and cross terms separately and consecutively. We write the signal term in the way
\begin{align}\label{decompdis}\frac{1}{i}\sum_{j=i}^n\d_j^{i}X(\d_j^{i}X)^{\top}&=\sum_{j=1}^n\d_jX(\d_jX)^{\top}+\sum_{l=1}^n\d_lX\sum_{j=1}^{i\wedge l}\left(1-\frac{j}{i}\right)(\d_{l-j}X)^{\top}\\
&\nonumber\quad +\sum_{l=1}^n\sum_{j=1}^{i\wedge l}\left(1-\frac{j}{i}\right)\d_{l-j}X(\d_lX)^{\top}\,.\end{align}
The first addend is the realized volatility matrix, converging with rate $n^{1/2}$ to the integrated volatility matrix, and thus contributing only an asymptotically negligible error term. Because of \[\vec\big(\d_lX(\d_{l-j}X)^{\top}+\d_{l-j}X(\d_{l}X)^{\top}\big)=\mathcal{Z} \vec\big(\d_lX(\d_{l-j}X)^{\top}\big)\,,\] it is enough to consider one addend. Using
\begin{align*}\cov\Big(\vec\big(\d_lX(\d_{l-j}X)^{\top}\big)\big|\mathcal{F}_{\frac{(l-1)T}{n}}\Big)&=\cov\Big(\d_{l-j} X\otimes \d_l X\big|\mathcal{F}_{\frac{(l-1)T}{n}}\Big)\\
&=\,\E\Big[\big(\d_{l-j} X\otimes \d_l X\big)\big((\d_{l-j} X)^{\top}\otimes (\d_l X)^{\top}\big)\big|\mathcal{F}_{\frac{(l-1)T}{n}}\Big]\\
&\asymp^p \d_{l-j} X(\d_{l-j} X)^{\top}\otimes \Sigma_{\frac{(l-1)T}{n}}\,\tfrac{T}{n}\,,\end{align*}
\begin{align*}\E\Big[\d_{l-j} X(\d_{l-j} X)^{\top}\otimes \Sigma_{\frac{(l-1)T}{n}}\,\tfrac{T}{n}\Big]=T^2n^{-2}\Big(\Sigma_{\frac{(l-j-1)T}{n}}\otimes \Sigma_{\frac{(l-1)T}{n}}\Big)\,,\end{align*}
and $\mathcal{Z}^2=2\mathcal{Z}$, we derive that
\begin{align*}&\sum_{l=1}^n\cov\Bigg(\sum_{j=1}^{i\wedge l}\left(1-\frac{j}{i}\right)\mathcal{Z}\vec\Big(\Delta_l X\big(\Delta_{l-j}X\big)^{\top}\Big)\big|\mathcal{F}_{\frac{(l-1)T}{n}}\Bigg)\\
&\quad \asymp^p \sum_{l=1}^nT^2n^{-2}\Big(\Sigma_{\frac{(l-1)T}{n}}\otimes \Sigma_{\frac{(l-1)T}{n}}\Big)\,2\mathcal{Z}\sum_{j=1}^{i\wedge l}\left(1-\frac{j}{i}\right)^2\,.
\end{align*}
The smoothness of $\Sigma$ ensured by Assumption \ref{eff} giving the bound in \eqref{prel1} suffices that the approximation errors by replacing $\d_{l-j} X(\d_{l-j} X)^{\top}$ with its expectation and by replacing $\Sigma_{\frac{(l-j-1)T}{n}}$ with $\Sigma_{\frac{(l-1)T}{n}}$ are asymptotically negligible.
From the above considerations and verifying all other conditions of Theorem 3--1 of \cite{jacod1}, which readily follow along the same lines as in the proof of Proposition A.3 of \cite{bibinger2}, we obtain first stable central limit theorems for discretization errors of subsampling estimators with fixed subsampling frequencies. The covariances between them are determined with
\begin{align*}&\sum_{l=1}^n\E\hspace*{-.05cm}\left[\sum_{j=1}^{i\wedge l}\hspace*{-.05cm}\left(1-\frac{j}{i}\right)\hspace*{-.05cm}\mathcal{Z}\hspace*{-.05cm}\vec\hspace*{-.025cm}\Big(\Delta_l X\big(\Delta_{l-j}X\big)^{\top}\Big)\Big(\sum_{j=1}^{i'\wedge l}\left(1-\frac{j}{i'}\right)\hspace*{-.05cm}\mathcal{Z}\hspace*{-.05cm}\vec\hspace*{-.025cm}\Big(\Delta_l X\big(\Delta_{l-j}X\big)^{\top}\Big)\Big)^{\top}\Big|\mathcal{F}_{\frac{(l-1)T}{n}}\hspace*{-.05cm}\right]\\
&\quad \asymp^p \sum_{l=1}^nT^2n^{-2}\Big(\Sigma_{\frac{(l-1)T}{n}}\otimes \Sigma_{\frac{(l-1)T}{n}}\Big)\,2\mathcal{Z}\sum_{j=1}^{\min{\left(i,i',l\right)}}\left(1-\frac{j}{i}\right)\left(1-\frac{j}{i'}\right)\,.
\end{align*}
Hence, we are left to evaluate the deterministic sum: \begin{align*}\sum_{j=1}^{m}\left(1-\frac{j}{i}\right)\left(1-\frac{j}{i'}\right)=\frac{m}{2}-\frac{m^2}{6\textsl{M}}-\frac{1}{8}+\frac{1}{12\textsl{M}}\asymp \frac{m}{6}\left(3-\frac{m}{\textsl{M}}\right)\,,\end{align*}
where $m=\min{\left(i,i'\right)}$ and $\textsl{M}=\max{\left(i,i'\right)}$. Including the weights according to \eqref{weights}, we set
\begin{align*}\mathfrak{D}^{\alpha}=\lim_{n\rightarrow\infty}{M_n^{-1}\sum_{k=1}^{M_n}\sum_{l=1}^k\frac{l}{6M_n}\left(3-\frac{l}{k}\right)\alpha_k\alpha_l}\,.\end{align*}
With the covariances for different subsample frequencies above we obtain a stable central limit theorem for vectors spanning over finite sets of different frequencies. The Cram\'{e}r-Wold device implies central limit theorems for linear combinations of the components. The final stable limit theorem
\begin{align}\label{hh1}n^{\frac14}\vec\Big(\sum_{i=1}^n\frac{\alpha_i}{i}\sum_{j=i}^n\Delta_j^iX\big(\Delta_j^i X\big)^{\top}-\int_0^T\Sigma_s\,ds\Big)\stackrel{st}{\rightarrow} MN\Big(0,4\mathfrak{D}^{\alpha}c T\int_0^T\big(\Sigma_s\otimes\Sigma_s\big)\mathcal{Z}\,ds\Big)\end{align}
is concluded by extending this to infinitely many subsample frequencies adopting the analogous step from \cite{zhang} for the univariate multi-scale estimator.
Thereby we conclude the signal term of \eqref{covms}. $\mathfrak{D}^{\alpha}$ is a constant showing up in the asymptotic discretization variance depending on the weights, where for the standard weights \eqref{cubweights} or cubic kernel $\mathfrak{D}^{\alpha}=13/70$. \\
Next, consider the noise term
\begin{align}\notag\sum_{i=1}^{M_n}\frac{\alpha_i}{i}\sum_{j=i}^n\left(\epsilon_j-\epsilon_{j-i}\right)\left(\epsilon_j^{\top}-\epsilon_{j-i}^{\top}\right)=\sum_{i=1}^{M_n}\frac{\alpha_i}{i}\left(2\sum_{j=1}^n\epsilon_j\epsilon_j^{\top}\right.\\
\left.-\sum_{j=i}^n\left(\epsilon_j\epsilon_{j-i}^{\top}+\epsilon_{j-i}\epsilon_{j}^{\top}\right)-\sum_{j=n-i+1}^n\epsilon_j\epsilon_j^{\top}-\sum_{j=0}^{i-1}\epsilon_j\epsilon_j^{\top}\right)\,.\label{noisedec}\end{align}
The last two sums lead for the non-adjusted multi-scale estimator \eqref{ms} to the negative bias by noise and end-effects.
The first inner sum on the right-hand side above does not depend on $i$ and the term vanishes since $\sum_{i=1}^{M_n}\alpha_i/i=0$. The variance-covariance matrices of the remaining uncorrelated addends contribute to the total variance-covariance matrix due to noise perturbation. As the noise variance-covariance matrix $\mathbf{H}$ is fixed, we may work conditional on $X$ and consider covariances directly instead of conditional covariances as for the discretization part. Denote the constant limits
\begin{align*}{\mathfrak{N}}^{\alpha}_2=\lim_{n\rightarrow\infty}{M_n\sum_{j=1}^{M_n-1}\left(\sum_{i=j+1}^{M_n}\frac{\alpha_i}{i}\right)^2} ~\mbox{and}~{\mathfrak{N}}^{\alpha}_1=\lim_{n\rightarrow\infty}{M_n^3\sum_{i=1}^{M_n}\frac{\alpha_i^2}{i^2}}\,.\end{align*}
\begin{align*}\text{Rewriting}\hspace*{3.15cm}\sum_{i=1}^{M_n}\frac{\alpha_i}{i}\sum_{j=i}^n\epsilon_j\epsilon_{j-i}^{\top}=\sum_{j=1}^n\sum_{i=1}^{M_n\wedge j}\frac{\alpha_i}{i}\epsilon_j\epsilon_{j-i}^{\top}\,,\hspace*{4.85cm}\end{align*}
\begin{align*}\sum_{i=1}^{M_n}\frac{\alpha_i}{i}\left(\sum_{j=0}^{i-1}\epsilon_j\epsilon_{j}^{\top}+\sum_{j=n-i+1}^{n}\epsilon_j\epsilon_{j}^{\top}\right)=\sum_{j=0}^{M_n-1}\left(\epsilon_j\epsilon_j^{\top}+\epsilon_{n-j}\epsilon_{n-j}^{\top}\right)\sum_{i=j+1}^{M_n}\frac{\alpha_i}{i}\,,\end{align*}
\begin{align}\label{h1}\text{we obtain that}\hspace*{1.375cm}\frac{M_n^3}{n}\cov\left(\vec\Bigg(\sum_{i=1}^{M_n}\frac{\alpha_i}{i}\sum_{j=i}^n\epsilon_j\epsilon_{j-i}^{\top}\Bigg)\right)\rightarrow \mathfrak{N}^{\alpha}_1\big(\bf{H}\otimes\bf{H}\big)\,,\hspace*{1.175cm}\end{align}
\begin{align}\label{h2}M_n\cov\left(\vec\Bigg(\sum_{i=1}^{M_n}\frac{\alpha_i}{i}\Bigg(\sum_{j=n-i+1}^n\epsilon_j\epsilon_j^{\top}+\sum_{j=0}^{i-1}\epsilon_j\epsilon_j^{\top}\Bigg)\Bigg)\right)\rightarrow  2 {\mathfrak{N}}^{\alpha}_2\,\mathbf{H}^{\otimes}\mathcal{Z}\,.\end{align}
Using once again that $\vec(A+A^{\top})=\mathcal{Z}\vec(A)$ for $A\in\R^{d\times d}$, we conclude the noise parts in \eqref{covms}. For the specific weights \eqref{cubweights} corresponding to the cubic kernel, we have $\mathfrak{N}^{\alpha}_2=6/5$ and $\mathfrak{N}^{\alpha}_1=12$, which gives the minimum of the variance due to noise, cf.\,\cite{zhang}.\\
Finally, consider the cross terms. They can be decomposed in addends of the form
\begin{align}\label{h4}\sum_{i=1}^{M_n}\frac{\alpha_i}{i}\sum_{j=0}^n\big(\zeta_{i,j}\epsilon_j^{\top}+\epsilon_j\zeta_{i,j}^{\top}\big)\,,\end{align}
\begin{align*}\mbox{where}~~\zeta_{i,j}=\begin{cases}-\d_{i-j}^i X&,0\le j\le (i-1)\\
\d_j^iX-\d_{j+i}^iX&,i\le j\le (n-i)\\
\d_j^i X &,n-i+1\le j\le n\end{cases}~.\end{align*}
In order to derive the asymptotic variance-covariance matrix, observe that
\begin{align*}\cov\Big(\vec\big(\zeta_{i,j}\epsilon_j^{\top}\big)\Big)&=\cov\Big(\epsilon_j\otimes \zeta_{i,j}\Big)=\E\Big[\big(\epsilon_j\otimes \zeta_{i,j}\big)\big(\epsilon_j^{\top}\otimes \zeta_{i,j}^{\top}\big)\Big]={\bf{H}}\otimes \E\big[\zeta_{i,j}\zeta_{i,j}^{\top}\big]\,.\end{align*}
Now, if we assume without loss of generality $1\le i\le i'\le M_n$, it holds for $M_n\le j\le (n-M_n)$ that
\begin{align*}\zeta_{i,j}\zeta_{i',j}^{\top}=\Bigg(\sum_{r=j-i}^{j-1}\d_rX-\sum_{r=j}^{i+j-1}\d_rX\Bigg)\Bigg(\sum_{r=j-i'}^{j-1}(\d_rX)^{\top}-\sum_{r=j}^{i'+j-1}(\d_rX)^{\top}\Bigg)\\
=\zeta_{i,j}\zeta_{i,j}^{\top}+\zeta_{i,j}\Bigg(\sum_{r=j-i'}^{j-i-1}(\d_rX)^{\top}-\sum_{r=j+i}^{i'+j-1}(\d_r X)^{\top}\Bigg)\,.\end{align*}
We obtain for the sum of conditional variance-covariance matrices the following convergence:
\begin{align}
&\notag M_n\sum_{j=0}^n\cov\left(\vec\Bigg(\sum_{i=1}^{M_n}\frac{\alpha_i}{i}\big(\zeta_{i,j}\epsilon_j^{\top}+\epsilon_j\zeta_{i,j}^{\top}\big)\Bigg)\Big|\mathcal{F}_{\frac{(j-1)T}{n}}\right)\\
&\notag\quad =M_n\sum_{j=0}^n\cov\left(\sum_{i=1}^{M_n}\frac{\alpha_i}{i}\mathcal{Z}\vec\big(\zeta_{i,j}\epsilon_j^{\top}\big)\Big|\mathcal{F}_{\frac{(j-1)T}{n}}\right)\\
 &\notag\quad\asymp^p 2M_n\mathcal{Z} \Bigg({\bf{H}}\otimes  \sum_{i=1}^{M_n}\sum_{r=1}^{M_n}\frac{\alpha_i\alpha_r}{ir}(i\wedge r)\Big(\frac{1}{(i\wedge r)}\sum_{j=(i\wedge r)}^n\E\Big[\d_j^{(i\wedge r)}X(\d_j^{(i\wedge r)}X)^{\top}\Big)\Big|\mathcal{F}_{\frac{(j-1)T}{n}}\Big]\Bigg) \mathcal{Z}\\
&\label{hh2}\quad \stackrel{p}{\rightarrow} 2\mathfrak{M}^{\alpha}\mathcal{Z}\Big({\bf{H}}\otimes \left[X,X\right]_T\Big)\mathcal{Z}\,,\end{align}
with the constant of the limit dependening on the weights \eqref{weights}:
\begin{align}\label{limitm}\mathfrak{M}^{\alpha}=\lim_{M_n\rightarrow\infty}\sum_{i,k=1}^{M_n}\frac{\alpha_i\alpha_k}{ik} (i\wedge k)\,.\end{align}
For the specific weights \eqref{cubweights}, the constant takes the value $\mathfrak{M}^{\alpha}=6/5$. To elucidate the structure of the asymptotic variance-covariance matrix of the cross terms we use the decomposition $\mathcal{Z}=E_{d^2}+C_{d,d}\in \mathds{R}^{d^2\times d^2}$ with $E_{d^2}$ the identity matrix and $C_{d,d}$ the so-called commutation matrix characterized by $C_{d,d}\vec(A)=\vec(A^{\top})$ for $A\in\mathds{R}^{d\times d}$ and satisfying $C_{d,d}(A\otimes B)C_{d,d}=(B\otimes A)$ for $A,B\in\mathds{R}^{d\times d}$. Then, we can show the following identity for $A,B\in\mathds{R}^{d\times d}$:
\begin{align}\label{kron}\mathcal{Z}(A\otimes B)\mathcal{Z}&=\big(E_{d^2}+C_{d,d}\big)(A\otimes B)\big(E_{d^2}+C_{d,d}\big)\\
&\notag=(A\otimes B)+(A\otimes B)C_{d,d}+C_{d,d}(A\otimes B)+C_{d,d}(A\otimes B)C_{d,d}\\
&\notag=(A\otimes B)+(A\otimes B)C_{d,d}+(B\otimes A)C_{d,d}+(B\otimes A)\\
&\notag=\big(A\otimes B+B\otimes A\big)\mathcal{Z}\\
&\notag=\mathcal{Z}(B\otimes A)\mathcal{Z}\,.\end{align}
This illuminates the different illustrations of the variance-covariance matrix of cross terms and that we can work equivalently with the transpose term above.\\
The remaining elements of the proof of a multivariate stable central limit theorem are close to \cite{zhang} and \cite{bibinger2} again founded on the multivariate stable convergence theorem by \cite{jacod1}. Hence, we restrict ourselves to the evaluation of the general multivariate variance-covariance structure and derive \eqref{covms} by \eqref{hh1}, \eqref{h1}, \eqref{h2} and \eqref{hh2}. This completes the proof of Theorem \ref{procovms}.$\hfill\Box$

\section{Proofs of Section \ref{sec:3}}
\subsection{Proof of Theorem \ref{procovgms} and Corollary \ref{endocor}}
The strategy of proof follows the one of Theorem \ref{procovms}, but generalizing the variance-covariance structure to irregular observation times. The proof of the stable central limit theorem is again traced back to Theorem 3--1 of \cite{jacod1}.\\
Let us begin with the discretization error which is the key step to integrate endogenous observation times. We condition on  $\big(\mathcal{F}_{T_{l-1}}\Big)_l$ below which is each time replaced by $\big(\mathcal{F}^N_{T_{l-1}} \big)_l$ defined within Assumption \ref{disendo} for the endogenous case.  Based on a decomposition analogously to \eqref{decompdis}, we find that
\begin{align*}&\frac{N}{M_n}\sum_{j=1}^N\cov\Big(\vec\Big(\sum_{i=1}^{M_N\wedge j}\frac{\alpha_i}{i}\big(X_{T_j}^+-X_{T_{j-i}}^-\big)\big(X_{T_j}^+-X_{T_{j-i}}^-\big)^{\top}\Big)\Big|\mathcal{F}_{T_{j-1}}\Big)\\
&\quad\asymp^p\frac{N}{M_n}\sum_{l=1}^N\cov\Bigg(\sum_{i=1}^{M_N\wedge l}{\alpha_i}\mathcal{Z}\vec\Big(\big(X_{T_l}-X_{T_{l-1}}\big)\sum_{j=1}^{i\wedge l}\left(1-\frac{j}{i}\right)\big(X_{T_{l-j}}-X_{T_{l-j-1}}\big)^{\top}\Big)\Big|\mathcal{F}_{T_{l-1}}\Bigg)\end{align*}
using $\vec(A+A^{\top})=\mathcal{Z}\vec(A)$ for $A\in\mathds{R}^{d\times d}$, that the realized volatility matrix has an asymptotically negligible discretization error as in \eqref{decompdis} and that next- and previous-tick interpolations are asymptotically negligible in the discretization variance-covariance matrix on Assumption \ref{dis4} or Assumption \ref{disendo}, respectively. The bounds
\begin{align*}X_{T_j}^+-X_{T_{j-i}}^-=X_{T_j}-X_{T_{j-i}}+\mathcal{O}_{\P}\Big(\max_{p,p'}\{(t_p^+(T_j)-T_j)^{1/2},(T_{j-i}-t_{p'}^-(T_{j-i}))^{1/2}\}\Big)\end{align*}
with \eqref{eqdis3} (resp.\,\eqref{eqdis4}) suffice to prove that interpolation terms do not trigger the asymptotic variance-covariance, see Proposition A.\,10 of \cite{bibinger2} for a rigorous proof that directly carries over.\\
As above, we can neglect the drift part of $X$
, such that
\begin{align}\notag&\cov\Big(\vec\Big((X_{T_l}-X_{T_{l-1}})(X_{T_{l-j}}-X_{T_{l-j-1}})^{\top}\Big)\Big)\\
&\notag\asymp \E\left[\int_{T_{l-j-1}}^{T_{l-j}}\sigma_{s}\,dW_s\Big(\int_{T_{l-j-1}}^{T_{l-j}}\sigma_{s}\,dW_s\Big)^{\top}\right]\otimes \E\left[\int_{T_{l-1}}^{T_{l}}\sigma_{s}\,dW_s\Big(\int_{T_{l-1}}^{T_{l}}\sigma_{s}\,dW_s\Big)^{\top}\right]\\
&\label{keyendo}\asymp \E\left[\int_{T_{l-j-1}}^{T_{l-j}}\Sigma_s\,ds\right]\otimes \E\left[\int_{T_{l-1}}^{T_{l}}\Sigma_s\,ds\right]\,,\end{align}
\vspace*{-.75cm}

\allowdisplaybreaks[4]{
\begin{align*}
\E\Bigg[\int_{T_{l-j-1}}^{T_{l-j}}\Sigma_s\,ds\otimes \int_{T_{l-1}}^{T_{l}}\Sigma_s\,ds\Big|\mathcal{F}_{T_{l-1}}\Bigg]\asymp^p \big(\Sigma_{T_{l-j-1}}\otimes \Sigma_{T_{l-1}}\big)(T_l-T_{l-1})(T_{l-j}-T_{l-j-1})\,,\end{align*}
by generalized It\^{o}-isometry. In particular, we only consider these terms for $j>0$ here and the above transformation applies to both, exogenous observation times on Assumption \ref{dis4} and endogenous observation times on Assumption \ref{disendo}. The proof for the discretization error follows the same strategy as for regular observation times above. We combine several steps in the following computation which is analogous for both, exogenous and endogenous observation times. As in the proof of Theorem \ref{procovms}, exploiting the smoothness of $\Sigma$ yields
\begin{align*}&\frac{N}{M_N}\sum_{j=1}^N\cov\Big(\vec\Big(\sum_{i=1}^{M_N\wedge j}\frac{\alpha_i}{i}\big(X_{T_j}^+-X_{T_{j-i}}^-\big)\big(X_{T_j}^+-X_{T_{j-i}}^-\big)^{\top}\Big)\Big|\mathcal{F}_{T_{j-1}}\Big)\\
&\quad\asymp^p\frac{N}{M_N}\sum_{l=1}^N\cov\Big(\mathcal{Z}\sum_{i=1}^{M_N\wedge l}\frac{\alpha_i}{i}\vec\Big((X_{T_l}-X_{T_{l-1}})\sum_{q=1}^{i\wedge l}\left(1-\frac{q}{i}\right)(X_{T_{l-q+1}}-X_{T_{l-q}})^{\top}\Big)\Big|\mathcal{F}_{T_{l-1}}\Big)\\
&\quad\asymp^p2\frac{N}{M_N}\sum_{l=1}^N\sum_{i,k=1}^{M_N\wedge l}\alpha_i\alpha_k(T_l-T_{l-1})\mathcal{Z}\Bigg(\sum_{q=1}^{\min{(l,i,k)}}\left(1-\frac{q}{i}\right)\left(1-\frac{q}{k}\right)\\
&\hspace*{6.5cm}\times (X_{T_{l-q+1}}-X_{T_{l-q}})(X_{T_{l-q+1}}-X_{T_{l-q}})^{\top}\otimes \Sigma_{T_{l-1}}\Bigg)\mathcal{Z}\\
&\quad\asymp^p2\frac{N}{M_N}\sum_{i,k=1}^{M_N}\alpha_i\alpha_k\sum_{l=1}^N(T_l-T_{l-1})\mathcal{Z}\big(\Sigma_{T_{l-1}}\otimes \Sigma_{T_{l-1}}\big)\mathcal{Z}\sum_{q=1}^{\min{(l,i,k)}}\left(1-\frac{q}{i}\right)\left(1-\frac{q}{k}\right)\,.\end{align*}
Inserting $A=B$ in \eqref{kron} above yields for $A\in\mathds{R}^{d\times d}$ the identity
\begin{align}\label{kron2}\mathcal{Z}(A\otimes A)\mathcal{Z}=2\mathcal{Z}(A\otimes A)=2(A\otimes A)\mathcal{Z}\,.\end{align}
Assumption \ref{aqvt2} (resp.\,Assumption \ref{disendo}) ensures the convergence
\begin{align*}&\frac{N}{M_N}\sum_{j=1}^N\cov\Big(\vec\Big(\sum_{i=1}^{M_N\wedge j}\frac{\alpha_i}{i}\big(X_{T_j}^+-X_{T_{j-i}}^-\big)\big(X_{T_j}^+-X_{T_{j-i}}^-\big)^{\top}\Big)\Big|\mathcal{F}_{T_{j-1}}\Big)\\
&\hspace*{6.5cm}\stackrel{p}{\rightarrow} 4\int_0^T (D^{\alpha})^{\prime}(s)\big(\Sigma_s\otimes \Sigma_s\big)\mathcal{Z}\,ds\,.\end{align*}
The existence of $(D^{\alpha})^{\prime}(s)$ is ensured on Assumption \ref{aqvt2} (resp.\,\ref{disendo}) by dominated convergence and the convergence of $M_N^{-1}\sum\alpha_i\alpha_k(i\wedge k)$ on the conditions for the weights \eqref{weights}.\\
In the remainder of the proof possible endogeneity of observation times plays a minor role. The variance-covariance terms hinge on certain characteristics of observation times for which we assume (stochastic) convergence by Assumption \ref{aqvt2} and Assumption \ref{disendo}, respectively.
Following an analogous decomposition of the error due to noise as above in \eqref{noisedec}, we obtain two terms generalizing \eqref{h1} and \eqref{h2}, respectively. Observe that
\begin{align*}\cov\Big(\vec\Big(\epsilon_{T_j}^+\big(\epsilon_{T_{j-i}}^-\big)^{\top}\Big)\Big)&=\E\left[\epsilon_{T_{j-i}}^-\big(\epsilon_{T_{j-i}}^-\big)^{\top}\right]
\otimes\E\left[\epsilon_{T_j}^+\big(\epsilon_{T_j}^+\big)^{\top}\right]\,,\end{align*}
and that $\E\big[\epsilon_{T_j}^+\big(\epsilon_{T_{j'}}^+\big)^{\top}\big]$ vanishes whenever $|j-j'|>1$ and $\E\big[\epsilon_{T_l}^-\big(\epsilon_{T_{l'}}^-\big)^{\top}\big]$ for all $l\ne l'$. In fact, successive next-ticks can coincide $t_p^+(T_j)=t_p^+(T_{j-1})$, while $t_p^-(T_j)=t_p^-(T_{j'})$ implies $j=j'$, see Section 4 of \cite{bibinger2} for a discussion of this aspect. Thereby, we deduce that
\begin{align*}\frac{M_N^3}{N}\cov\Big(\vec\Big(\sum_{i=1}^{M_N}\frac{\alpha_i}{i}\sum_{j=i}^N\epsilon_{T_j}^+\big(\epsilon_{T_{j-i}}^-\big)^{\top}\Big)\Big)\rightarrow \big(\mathbf{H}\otimes \mathbf{H}\big)^*\,,\end{align*}
as generalization of \eqref{h1} with \eqref{starlimit} as well as a simple upper bound for $|\alpha_i^2/i^2-\alpha_i\alpha_{i+1}/(i(i+1))|$.
The convergence assumption \eqref{2starlimit} readily gives the generalization of \eqref{h2}:
\begin{align*}M_N\cov\Bigg(\vec\Big(\sum_{i=1}^{M_N}\frac{\alpha_i}{i}\Big(\sum_{j=0}^{i-1}\epsilon_{T_j}^+\big(\epsilon_{T_j}^+\big)^{\top}+\sum_{j=N-i+1}^{N}\epsilon_{T_j}^-\big(\epsilon_{T_j}^-\big)^{\top}\Big)\Big)\Bigg)\rightarrow\mathfrak{N}_2 \,\tilde{\mathbf{H}}^{\otimes}\mathcal{Z}\,,\end{align*}
together with a simple upper bound on $\alpha_i^2/i^2-\alpha_i\alpha_{i+1}/(i(i+1))$.
We are left to consider the cross terms
\begin{align*}\sum_{i=1}^{M_N}\frac{\alpha_i}{i}\sum_{j=i}^N\Big(\Big(\epsilon_{T_j}^+-\epsilon_{T_{j-i}}^-\Big)\Big(X_{T_j}^+-X_{T_{j-i}}^-\Big)^{\top}+\Big(X_{T_j}^+-X_{T_{j-i}}^-\Big)\Big(\epsilon_{T_j}^+-\epsilon_{T_{j-i}}^-\Big)^{\top}\Big)\,.\end{align*}
The (only) main difference to the regular observation setup is that for any subsample-lag $i$ above all $\epsilon_j,i\le j\le N-i,$ occurred once as left and once as right end point of sub-sampled intervals, which is not necessarily the case here. We differentiate all times for which $t_p^+(T_j)=t_{p'}^+(T_{j-1})$ and $t_p^+(T_j)=t_{p'}^-(T_{u})$ for some $u$ and $p,p'\in\{1,\ldots,d\}$. Again, we use that successive next-ticks can coincide for several components (but never for all at the same time) while previous-ticks change each time. We derive that
\begin{align*}&M_N\sum_{j=1}^N\cov\Bigg(\mathcal{Z}\vec\Big(\sum_{i=1}^{M_N\wedge j}\frac{\alpha_i}{i}\Big(\epsilon_{T_j}^+-\epsilon_{T_{j-i}}^-\Big)\Big(X_{T_j}^+-X_{T_{j-i}}^-\Big)^{\top}\Big)\Big|\mathcal{F}_{T_{j-1}}\Bigg)\\
&\quad\asymp^p M_N\sum_{i,r=1}^{M_N}\frac{\alpha_i\alpha_r}{ir}\sum_{j=i\wedge r}^N\frac{i\wedge r}{N}\mathcal{Z}\Big(\Big(N\big(S^N(T_j)-S^N(T_{j-1})\big)^{(pq)}\mathbf{H}^{(pq)}\Big)_{1\le p,q\le d}\otimes \Sigma_{T_{j-1}}\Big)\mathcal{Z}\\
&\quad\asymp^p M_N\sum_{i,r=1}^{M_N}\frac{\alpha_i\alpha_r}{ir}(i\wedge r)\sum_{j=i\wedge r}^N\mathcal{Z}\Big(\Big(\big(S^{\prime}(T_{j-1})\big)^{(pq)}\mathbf{H}^{(pq)}\Big)_{1\le p,q\le d}\otimes \Sigma_{T_{j-1}}\Big)\mathcal{Z}(T_j-T_{j-1})\,,\end{align*}
such that with \eqref{limitm}, \eqref{kron} and \eqref{tildelimit} on Assumption \ref{aqvt2} (resp.\,\ref{disendo}) we conclude \eqref{covgms}. The remaining elements of the proof of a multivariate stable central limit follow similar as in \cite{bibinger2} founded on the multivariate stable convergence theorem by \cite{jacod1} and we omit them here. This completes the proof of Theorem \ref{procovgms} and Corollary \ref{endocor}.$\hfill\Box$}

\subsection*{Proof of Proposition \ref{covest}}
Proposition \ref{covest} follows in the same way as Proposition 5.1 in \cite{bibinger2}, extending terms to the multivariate notion given in \eqref{covestgms}. Consistency of the estimator $\hat{\mathbf{H}}$ is easily proved, \cite{zhangmykland} provide a central limit theorem with $\sqrt{n}$-convergence rate in the univariate setting. Under convergence assumptions \eqref{lasa} and \eqref{dra} the series with available observation times inserted will converge accordingly. Then, consistency of the local binwise multi-scale estimators yields consistency of the overall estimator.$\hfill\Box$

\bibliographystyle{mychicago}
\bibliography{literatur}

\vspace*{1cm}\noindent
Markus Bibinger,
Institut f\"ur Mathematik, Humboldt-Universit\"at zu Berlin,\\
Unter den Linden 6, 10099 Berlin, Germany\\
bibinger@math.hu-berlin.de\\[.5cm]
Per A. Mykland, Department of Statistics, The University of Chicago,\\
5734 University Avenue, Chicago, Illinois 60637, USA\\
mykland@pascal.uchicago.edu

\end{document}